# HJB EQUATIONS FOR CERTAIN SINGULARLY CONTROLLED DIFFUSIONS

By Rami Atar[1], Amarjit Budhiraja[2] and Ruth J. Williams[3]

*Technion–Israel Institute of Technology, University of North Carolina and University of California*

Given a closed, bounded convex set $\mathcal{W} \subset \mathbb{R}^d$ with nonempty interior, we consider a control problem in which the state process $W$ and the control process $U$ satisfy

$$W_t = w_0 + \int_0^t \vartheta(W_s)\,ds + \int_0^t \sigma(W_s)\,dZ_s + GU_t \in \mathcal{W}, \qquad t \ge 0,$$

where $Z$ is a standard, multi-dimensional Brownian motion, $\vartheta, \sigma \in C^{0,1}(\mathcal{W})$, $G$ is a fixed matrix, and $w_0 \in \mathcal{W}$. The process $U$ is locally of bounded variation and has increments in a given closed convex cone $\mathcal{U} \subset \mathbb{R}^p$. Given $g \in C(\mathcal{W})$, $\kappa \in \mathbb{R}^p$, and $\alpha > 0$, consider the objective that is to minimize the cost

$$J(w_0, U) \doteq \mathbb{E}\left[\int_0^\infty e^{-\alpha s} g(W_s)\,ds + \int_{[0,\infty)} e^{-\alpha s}\,d(\kappa \cdot U_s)\right]$$

over the admissible controls $U$. Both $g$ and $\kappa \cdot u$ ($u \in \mathcal{U}$) may take positive and negative values. This paper studies the corresponding dynamic programming equation (DPE), a second-order degenerate elliptic partial differential equation of HJB-type with a state constraint boundary condition. Under the controllability condition $G\mathcal{U} = \mathbb{R}^d$ and the finiteness of $\mathcal{H}(q) = \sup_{u \in \mathcal{U}_1}\{-Gu \cdot q - \kappa \cdot u\}$, $q \in \mathbb{R}^d$, where $\mathcal{U}_1 = \{u \in \mathcal{U} : |Gu| = 1\}$, we show that the cost, that involves an improper integral, is well defined. We establish the following: (i) the value function for the control problem satisfies the DPE (in the viscosity sense), and (ii) the condition $\inf_{q \in \mathbb{R}^d} \mathcal{H}(q) < 0$ is necessary and sufficient for uniqueness of solutions to the DPE. The existence and

Received November 2006; revised April 2007.
[1]Supported in part by Israel Science Foundation Grant 126/02 and NSF Grant DMS-06-00206.
[2]Supported in part by ARO Grants W911NF-04-1-0230 and W911NF-0-1-0080.
[3]Supported in part by NSF Grant DMS-06-04537.

uniqueness of solutions are shown to be connected to an intuitive "no arbitrage" condition.

Our results apply to Brownian control problems that represent formal diffusion approximations to control problems associated with stochastic processing networks.

**1. Introduction.** This paper studies a class of singular control problems for diffusions with state constraints. Such a problem involves finite dimensional processes denoted by $W$ and $U$, referred to as *state* and *control* processes, respectively. They satisfy

$$W_t = w_0 + \int_0^t \vartheta(W_s)\,ds + \int_0^t \sigma(W_s)\,dZ_s + GU_t, \qquad t \ge 0, \tag{1.1}$$

where $W$ and $U$ are adapted to a filtration for which $Z$ is a standard multidimensional Brownian motion, $\vartheta$ and $\sigma$ are Lipschitz continuous functions and $G$ is a fixed $d \times p$ matrix ($d \le p$). The sample paths of the control $U$ are RCLL (right continuous with finite left limits), and they are locally of bounded variation. In addition, these sample paths have increments in a given closed convex cone $\mathcal{U} \subset \mathbb{R}^p$. By "state constraints" we mean that for a control to be regarded as *admissible*, it is required that the corresponding state process stays for all time within a given closed, bounded convex set $\mathcal{W} \subset \mathbb{R}^d$ that has a nonempty interior. For existence of admissible controls, we need the following controllability condition which is assumed throughout:

$$\text{the set } G\mathcal{U} \doteq \{Gu : u \in \mathcal{U}\} \text{ equals } \mathbb{R}^d. \tag{1.2}$$

The control problem consists of minimizing a cost of the form

$$J(w_0, U) \doteq \mathbb{E}\left[\int_0^\infty e^{-\alpha s} g(W_s)\,ds + \int_{[0,\infty)} e^{-\alpha s}\,d(\kappa \cdot U_s)\right] \tag{1.3}$$

over all admissible controls $U$. Here, $g$ is a continuous function, $\kappa$ is a fixed vector in $\mathbb{R}^p$ and $\alpha > 0$. One of the main new aspects in this work is that $\kappa \cdot u$ may take both positive and negative values for different $u \in \mathcal{U}$, hence, the cost definition involves an improper integral (see [1] for a related treatment involving only proper integrals and for additional references on singular control of diffusions with state constraints). Define

$$\mathcal{H}(q) = \sup_{u \in \mathcal{U}_1} \{-Gu \cdot q - \kappa \cdot u\}, \qquad q \in \mathbb{R}^d, \tag{1.4}$$

where $\mathcal{U}_1 = \{u \in \mathcal{U} : |Gu| = 1\}$ and $|\cdot|$ denotes the Euclidean norm on $\mathbb{R}^d$. We assume throughout that

$$\mathcal{H}(q) < \infty \qquad \text{for some } q \in \mathbb{R}^d, \tag{1.5}$$

(equivalently, for all $q \in \mathbb{R}^d$, see Lemma A.1 in the Appendix). This assumption ensures that the integrals and expectation in (1.3) are well defined, and that an integration by parts formula is valid for the last integral in (1.3), as shown in Lemma 2.1. Note that (1.2) implies that $\mathcal{U}_1$ is nonempty and, thus, using (1.2) and (1.5), we have that $\mathcal{H}(q) \in (-\infty, \infty)$ for all $q \in \mathbb{R}^d$. Our main goal is to study existence and uniqueness for the corresponding Dynamic Programming Equation (DPE)

$$((\mathcal{L} + \alpha)\psi - g) \vee \mathcal{H}(D\psi) = 0, \tag{1.6}$$

with a boundary condition stipulating that the (viscosity sense) solution is a supersolution up to the boundary (this is known as the state constraint boundary condition). Here, for a smooth function $f$ defined on $\mathcal{W}$, $Df$ denotes the gradient of $f$, $D^2 f$ the Hessian matrix of $f$, $\mathcal{L}f = -\frac{1}{2}\text{trace}(\sigma\sigma' D^2 f) - \vartheta \cdot Df$. Because of the way the boundary condition is set, it is seen to be redundant if the function $\mathcal{H}$ assumes only nonnegative values, and therefore one might expect the condition

$$\inf_{q \in \mathbb{R}^d} \mathcal{H}(q) < 0 \tag{1.7}$$

to be important for uniqueness of solutions. Our main result (Theorem 2.1) establishes that the value function of the control problem satisfies the DPE and that (1.7) is necessary and sufficient for uniqueness of viscosity solutions.

We then investigate the role played by the condition

$$\{u \in \mathcal{U} : Gu = 0 \text{ and } \kappa \cdot u \le 0\} = \{0\}, \tag{1.8}$$

that we refer to as a "no arbitrage" condition (following the terminology for an analogous condition, Assumption 2.2 of [7], that we restate here as Assumption 3.2). This condition (1.8) states that there is no nontrivial control that maintains the current state without incurring an immediate increase in the cost. Lemma A.1 shows that (1.5) can be written in the following equivalent form, which is a slightly weaker condition than (1.8) (see Lemma A.2):

$$\{u \in \mathcal{U} : |Gu| \le \varepsilon \text{ and } \kappa \cdot u \le -1\} = \varnothing \quad \text{for all } \varepsilon > 0 \text{ sufficiently small.} \tag{1.9}$$

In particular, under (1.9), there exist solutions to the DPE. We also show (in Theorem 2.2) that (1.8), together with a mild condition on the cone $\mathcal{U}$ (2.15), implies (1.7), thus providing a sufficient condition for uniqueness of solutions for the DPE. Example 2.1(a) shows a case where (1.9) (and in particular, existence of solutions) holds, while (1.8) (and uniqueness) fails.

A principal motivation for this paper comes from a family of models referred to in the literature as Brownian control problems (BCPs), that arise as formal diffusion approximations to control problems associated with stochastic processing networks. We refer the reader to [6] for a detailed account on relationships between stochastic processing networks and BCPs.

In [7] the BCPs introduced in [6] were shown to be equivalent to singular control problems of the form studied in the current paper. A key assumption for this equivalence is what we have stated here as Assumption 3.2. A simple argument shows that this assumption, in fact, implies (1.8) (cf. Lemma 3.1). Consequently, our results provide a characterization of the value function associated with the singular control problem of [7] as the unique solution of a corresponding DPE (cf. Theorem 3.1).

We now comment on some aspects of the technique. The proof that the value function of the control problem is a solution of the DPE (1.6) with a state constraint boundary condition is obtained by establishing a dynamic programming principle (DPP) (see Proposition 4.1). For a class of singular control problems with state constraint, a similar DPP was recently obtained in [1]. The proof there crucially used certain monotonicity properties of the singular control term in the cost (see the proof of Lemma 8.1 of [1]) that are not valid in the current formulation. To treat the more general form of the cost function considered in the current paper, we give a different (and more direct) proof that does not appeal to monotonicity requirements.

It is well understood that an appropriate framework for second order degenerate elliptic equations, of which the DPE (1.6) is a special case, is through the theory of viscosity solutions. The paper [3] gives an excellent tutorial on the subject. However, typical comparison results in this theory which are used to argue uniqueness of solutions rely on a key coercivity property which is usually unavailable for DPEs corresponding to singular control problems. More precisely, writing the equation (1.6) in the form $\check{F}(x,\psi,D\psi,D^2\psi)=0,\ x\in\mathcal{W}$, for a suitable $\check{F}:\mathcal{W}\times\mathbb{R}\times\mathbb{R}^d\times\mathcal{S}(d)\to\mathbb{R}$, where $\mathcal{S}(d)$ denotes the space of symmetric $d\times d$ matrices, the standard coercivity condition (see [3], equation (3.13)) requires, for some $\gamma\in(0,\infty)$,

$$\check{F}(x,r,q,X)-\check{F}(x,s,q,X)\ge\gamma(r-s)\qquad \text{for } r\ge s,\ (x,q,X)\in\mathcal{W}\times\mathbb{R}^d\times\mathcal{S}(d).$$

This condition is clearly not satisfied, in general, for $\check{F}$ as in (1.6). As suggested in Section 5C of [3], the existence of a strict subsolution to (1.6) (i.e., a subsolution defined with a strict inequality, uniform over $\mathcal{W}$) may be used as a substitute for coercivity in the comparison argument. As noted earlier, condition (1.7) is necessary and sufficient for uniqueness. The role it plays in the comparison proof is precisely by enabling the construction of a strict subsolution to (1.6).

The proof that (1.7) follows from the no arbitrage condition (1.8) and (2.15) is surprisingly indirect (see Theorem 2.2). Although these conditions are purely algebraic (involving only $G$, $\kappa$ and $\mathcal{U}$), we have not found a direct proof. Our proof, in fact, relies on regularity of the value function of the control problem (1.1)–(1.3) (and perturbations thereof).

The paper is organized as follows. In Section 2 we present the basic setup and the main result (Theorem 2.1) showing existence of solutions for the DPE and giving necessary and sufficient conditions for uniqueness of solutions. We also establish the finiteness and Lipschitz continuity of the value function, and show how the hypothesis for uniqueness can be verified by means of the no arbitrage condition (1.8) and (2.15). Section 3 describes the connection of the control problem formulated in Section 2 with the reduced Brownian control problem identified in [7], and characterizes the value function of the latter in terms of the DPE. The latter is done by verifying the conditions of the main theorem. Section 4 is devoted to the proof of the DPP and establishing the solvability of the DPE by the value function. Finally, in Section 5 we establish uniqueness via a comparison result (Theorem 5.1), formulated for a class of equations that is, in fact, more general in terms of conditions on $\mathcal{W}$ and $\mathcal{H}$ than that of Section 2.

The following notation and terminology will be used.

We shall use $c$, $c_1, c_2, \dots$ to denote positive deterministic constants whose values may change from the proof of one result to another. For $\alpha, \beta \in \mathbb{R}^n$, $|\alpha|$ denotes the Euclidean norm of $\alpha$, and $\alpha \cdot \beta$ denotes the usual scalar product between $\alpha$ and $\beta$. The operator norm of a matrix $M$ will be denoted by $|M|$. Let $B_\varepsilon(x) = \{y \in \mathbb{R}^n : |x-y| < \varepsilon\}$ and let $S^{n-1}$ denote the unit sphere in $\mathbb{R}^n$. For a set $A \subset \mathbb{R}^n$, $A^o$ [$\overline{A}$, $\partial A$] denotes the interior [resp., closure, boundary] of $A$. The infimum over an empty set is regarded as $+\infty$. For a set $S \subset \mathbb{R}^n$, $C(S)$ [$C^2(S)$] denotes the space of continuous [resp. twice continuously differentiable] real valued functions defined on $S$.

For a function $f : [0,\infty) \to \mathbb{R}^n$ and $t \geq 0$, we write $|f|_t$ for the total variation of $f$ over $[0,t]$ with respect to the Euclidean norm, defined by

$$|f|_t = |f(0)| + \sup\left\{\sum_{i=1}^{l} |f(t_i) - f(t_{i-1})| : 0 = t_0 < t_1 < \cdots < t_l = t,\, l \geq 1\right\},$$

and $|f|_t^* = \sup_{s\in[0,t]} |f(s)|$. A function from $[0,\infty)$ to some Polish space $E$ is RCLL if it is right-continuous on $[0,\infty)$ and has left limits in $E$ on $(0,\infty)$. When $E \subset \mathbb{R}^n$, for such an RCLL $\xi$ denote $\Delta\xi(t) = \xi(t) - \xi(t-)$ for $t > 0$. As a convention, we set $\Delta\xi(0) \doteq \xi(0)$. The space of all RCLL functions from $[0,\infty)$ into $E$ is denoted by $\mathcal{D}([0,\infty), E)$ and is endowed with the usual Skorohod topology. The space of all continuous functions from $[0,\infty)$ into $E$ is denoted by $\mathcal{C}([0,\infty), E)$ and is endowed with the topology of uniform convergence on compacts. A function $u : [0,\infty) \to \mathbb{R}^n$ is said to *have increments in a set* $\mathcal{U}$ if $u(0) \in \mathcal{U}$ and $u(t) - u(s) \in \mathcal{U}$ for all $0 \leq s < t < \infty$.

An *$n$-dimensional process* is a measurable map from a measurable space $(\Omega, \mathcal{G})$ to $\mathcal{D}([0,\infty), \mathbb{R}^n)$, and a *process* is an $n$-dimensional process for some $n$. A *continuous* process is a process having continuous sample paths almost surely. For a process $X$, we use the notation $X(t)$ and $X_t$ interchangeably.

A *filtered probability space* is a quadruple $(\Omega, \mathcal{G}, \{\mathcal{G}_t\}, \mathbb{Q})$, where $(\Omega, \mathcal{G}, \mathbb{Q})$ is a probability space and $\{\mathcal{G}_t\}$ is a *filtration*, that is, a family of sub-$\sigma$-algebras of the $\sigma$-algebra $\mathcal{G}$ indexed by $t \in \mathbb{R}_+$ and satisfying $\mathcal{G}_s \subset \mathcal{G}_t$ whenever $0 \leq s < t < \infty$. An $n$-dimensional process $X = \{X(t) : t \in \mathbb{R}_+\}$ defined on such a filtered probability space is said to be *adapted* if for each $t \geq 0$ the function $X(t) : \Omega \to \mathbb{R}^n$ is measurable when $\Omega$ has the $\sigma$-algebra $\mathcal{G}_t$ and $\mathbb{R}^n$ has its Borel $\sigma$-algebra.

Given a positive integer $n$ and a filtered probability space $(\Omega, \mathcal{G}, \{\mathcal{G}_t\}, \mathbb{Q})$, a process $Z$ defined on this space is said to be an *$n$-dimensional $\{\mathcal{G}_t\}$-standard Brownian motion* if it is a continuous, adapted $n$-dimensional process such that:

(i) $Z_0 = 0$, $\mathbb{Q}$-a.s., and
(ii) for $0 \leq s < t$, under $\mathbb{Q}$, the increment $Z_t - Z_s$ is independent of $\mathcal{G}_s$ and is normally distributed with mean zero and covariance matrix $(t-s)I$, where $I$ stands for the $n \times n$ identity matrix.

## 2. Setting and main result.

Let $k$ be a given positive integer. We say that $\Phi = (\Lambda, \mathcal{G}, \{\mathcal{G}_t\}, \mathbb{Q}, Z)$ is a *system* if $(\Lambda, \mathcal{G}, \{\mathcal{G}_t\}, \mathbb{Q})$ is a filtered probability space endowed with a $k$-dimensional $\{\mathcal{G}_t\}$-standard Brownian motion $Z$. We consider a control problem in which the *state* process is to remain in a given closed, bounded, convex set $\mathcal{W} \subset \mathbb{R}^d$ that has a nonempty interior, and the *control* process is to have increments in a given nonempty closed convex cone $\mathcal{U} \subset \mathbb{R}^p$, where $d$ and $p$ are given positive integers satisfying $p \geq d$. A $d \times p$ matrix $G$ specifies the linear effect of this control process on the state process. Conditions (1.2) and (1.5) [recall the definition of $\mathcal{H}$ given in (1.4)] are assumed throughout this paper. The drift and diffusion coefficients are denoted by $\vartheta$, a function from $\mathcal{W}$ to $\mathbb{R}^d$, and $\sigma$, a function from $\mathcal{W}$ to the space of $d \times k$ matrices, respectively. Throughout, $\vartheta$ and $\sigma$ are assumed to be Lipschitz continuous on $\mathcal{W}$.

DEFINITION 2.1. Given $k, d, p, \mathcal{W}, \mathcal{U}, G, \vartheta, \sigma$ as described above, an admissible control for the initial condition $w_0 \in \mathcal{W}$ is a $p$-dimensional adapted process $U$ defined on some filtered probability space $(\Omega, \mathcal{G}, \{\mathcal{G}_t\}, \mathbb{P})$ for which there exists a $d$-dimensional adapted process $W$ and a $k$-dimensional $\{\mathcal{G}_t\}$-standard Brownian motion $Z$, such that the following three properties hold $\mathbb{P}$-a.s.:

(i) One has

$$W_t = w_0 + \int_0^t \vartheta(W_s)\, ds + \int_0^t \sigma(W_s)\, dZ_s + GU_t, \qquad t \geq 0; \tag{2.1}$$

(ii) $U$ is locally of bounded variation and has increments in $\mathcal{U}$;
(iii) $W(t) \in \mathcal{W}$ for all $t \geq 0$.

We call $W$ the state process corresponding to the control $U$. The class of all admissible controls for the initial condition $w_0$ will be denoted by $\mathcal{A}(w_0)$.

Since $\mathcal{U}$ is a convex cone, we have $U_t \in \mathcal{U}$ for all $t \geq 0$, a.s. Also, it is easy to see that $\mathcal{A}(w)$ is nonempty for all $w \in \mathcal{W}$ (see proof of Proposition 2.1).

REMARK 2.1. Note that when we select $U \in \mathcal{A}(\omega_0)$, it is implicit that $U$ carries with it a filtered probability space $(\Omega, \mathcal{G}, \{\mathcal{G}_t\}, \mathbb{P})$ and processes $W, Z$. Expectations under $\mathbb{P}$ will be denoted by $\mathbb{E}$ and we shall often write a.s. instead of $\mathbb{P}$-a.s.

Let $\kappa \in \mathbb{R}^p$ and $\alpha \in (0, \infty)$ be given. The cost criterion will involve an improper integral in the form of the right-hand side of equation (2.3) below. Note first that for all $w_0 \in \mathcal{W}$, $U \in \mathcal{A}(w_0)$ and $t > 0$,

$$\alpha \int_0^t e^{-\alpha s} \kappa \cdot U_s \, ds + e^{-\alpha t} \kappa \cdot U_t = \int_{[0,t]} e^{-\alpha s} \, d(\kappa \cdot U_s) \qquad \text{a.s.} \tag{2.2}$$

Here we recall the convention that the contribution to the integral on the right-hand side above at time zero is $\kappa \cdot U_0$. This identity was previously noted in [7] (see Lemma A.1 there which follows from Theorem 18, page 278 and Theorem 8, page 265 of [5]). The following lemma ensures that our cost functional is bounded below and that integration by parts holds on $[0, \infty)$.

LEMMA 2.1. (i) *One has*

$$\sup_{w_0 \in \mathcal{W}} \sup_{U \in \mathcal{A}(w_0)} \mathbb{E}\left[\int_0^\infty e^{-\alpha s} (\kappa \cdot U_s)^- \, ds\right] < \infty.$$

(ii) *The limit* $\lim_{t \to \infty} \int_{[0,t]} e^{-\alpha s} \, d(\kappa \cdot U_s)$, *denoted as* $\int_{[0,\infty)} e^{-\alpha s} \, d(\kappa \cdot U_s)$, *exists as an improper integral a.s. with values in* $(-\infty, \infty]$, *for each* $w_0 \in \mathcal{W}$ *and* $U \in \mathcal{A}(w_0)$. *Furthermore, the following* "*integration by parts*" *formula holds:*

$$\alpha \int_0^\infty e^{-\alpha s} \kappa \cdot U_s \, ds = \int_{[0,\infty)} e^{-\alpha s} \, d(\kappa \cdot U_s) \qquad \textit{a.s.} \tag{2.3}$$

Note that the integral on the left-hand side of (2.3) is well defined a.s. (possibly taking the value $+\infty$) by the first part of the lemma.

PROOF. From (1.5) and Lemma A.1 we have that for some $c_1 \in (0, \infty)$, $(\kappa \cdot u)^- \leq c_1 |Gu|$ for all $u \in \mathcal{U}$. Thus, from (2.1) and compactness of $\mathcal{W}$, we have, for $U \in \mathcal{A}(w_0)$,

$$\mathbb{E}\left[\int_0^\infty e^{-\alpha s} (\kappa \cdot U_s)^- \, ds\right]$$

$$\begin{aligned}
&\le c_1\mathbb{E}\left[\int_0^\infty e^{-\alpha s}\left|W_s - w_0 - \int_0^s \vartheta(W_u)\,du - \int_0^s \sigma(W_u)\,dZ_u\right| ds\right] \qquad (2.4)\\
&\le c_1\left[\int_0^\infty e^{-\alpha s}(2c_2 + \|\vartheta\|_{\mathcal{W}} s + \|\sigma\|_{\mathcal{W}}\sqrt{s})\,ds\right],
\end{aligned}$$

where $c_2 = \sup\{|w| : w \in \mathcal{W}\}$, $\|\vartheta\|_{\mathcal{W}} = \sup\{|\vartheta(w)| : w \in \mathcal{W}\}$, $\|\sigma\|_{\mathcal{W}} = \sup\{|\sigma(w)| : w \in \mathcal{W}\}$ and we have used an $L^2$-isometry for stochastic integrals to obtain one of the bounds. Since the last integral above is finite and independent of $w_0 \in \mathcal{W}$ and $U \in \mathcal{A}(w_0)$, (i) holds. The proof of (ii) follows the argument given in Lemmas 3.2 and A.4 of [7] and for completeness is included in the Appendix. □

Let $g : \mathcal{W} \to \mathbb{R}$ be a continuous function. For $w_0 \in \mathcal{W}$ and $U \in \mathcal{A}(w_0)$, let the associated cost be defined as

$$J(w_0, U) \doteq \mathbb{E}\left[\int_0^\infty e^{-\alpha s} g(W_s)\,ds + \int_{[0,\infty)} e^{-\alpha s}\,d(\kappa\cdot U_s)\right]. \qquad (2.5)$$

In view of Lemma 2.1(ii), we have the following equivalent representation for the cost:

$$J(w_0, U) = \mathbb{E}\left[\int_0^\infty e^{-\alpha s} g(W_s)\,ds + \alpha\int_0^\infty e^{-\alpha s}\kappa\cdot U_s\,ds\right]. \qquad (2.6)$$

In view of the boundedness of $g$ on $\mathcal{W}$ and Lemma 2.1(i), $J(w_0, U)$ is bounded below uniformly for $w_0 \in \mathcal{W}$ and $U \in \mathcal{A}(w_0)$. The value function is defined as

$$V(w_0) \doteq \inf_{U\in\mathcal{A}(w_0)} J(w_0, U), \qquad w_0 \in \mathcal{W}. \qquad (2.7)$$

The following elementary lemma will be used at several places.

LEMMA 2.2. *There is a Borel measurable function $\varpi : \mathbb{R}^d \to \mathcal{U}$ and $c_\varpi \in (0,\infty)$ such that $G\varpi(x) = x$ and $|\varpi(x)| \le c_\varpi|x|$ for all $x \in \mathbb{R}^d$.*

PROOF. Let $\{e_i\}_{i=1}^d$ be an orthonormal basis in $\mathbb{R}^d$. From (1.2) there exist $f_i^+, f_i^- \in \mathcal{U}$ such that $Gf_i^+ = e_i$ and $Gf_i^- = -e_i$, $i = 1,\ldots,d$. Therefore,

$$\varpi(x) \doteq \sum_{i=1}^d |\langle x, e_i\rangle|(f_i^+ 1_{\{\langle x,e_i\rangle>0\}} + f_i^- 1_{\{\langle x,e_i\rangle\le 0\}}), \qquad x \in \mathbb{R}^d,$$

satisfies all of the desired properties with $c_\varpi \doteq d\max_{i=1}^d(|f_i^+| \vee |f_i^-|)$. □

Let $C^2(\mathcal{W})$ denote the set of twice continuously differentiable functions defined from $\mathcal{W}$ into $\mathbb{R}$. Let $\Gamma = \sigma\sigma'$ and, for $f \in C^2(\mathcal{W})$, let

$$\mathcal{L}f(x) = -\frac{1}{2}\mathrm{trace}(\Gamma(x)D^2 f(x)) - \vartheta(x)\cdot Df(x)$$

(2.8)
$$= -\frac{1}{2}\sum_{i,j=1}^{d}\Gamma_{ij}(x)\frac{\partial^2 f}{\partial x_i\,\partial x_j}(x) - \sum_{i=1}^{d}\vartheta_i(x)\frac{\partial f}{\partial x_i}(x), \qquad x\in\mathcal{W}.$$

The following proposition gives some basic properties of the value function.

PROPOSITION 2.1. *For each $w_0 \in \mathcal{W}$, the value is finite, that is, $V(w_0) \in (-\infty,\infty)$. The value function $V:\mathcal{W}\to\mathbb{R}$ is Lipschitz continuous, that is, there exists a constant $c_{\text{lip}}\in(0,\infty)$ such that, for all $w_1,w_2\in\mathcal{W}$, $|V(w_1)-V(w_2)|\le c_{\text{lip}}|w_1-w_2|$.*

PROOF. Note first that $V(w_0) > -\infty$ for $w_0\in\mathcal{W}$ is immediate from Lemma 2.1(i). To show that $V(w_0)<\infty$, consider first $w_0\in\mathcal{W}^o$. Let $B=\overline{B_\varepsilon(w_0)}$, where $\varepsilon$ is so small that $B\subset\mathcal{W}$. Let $n(w)$ denote the inward unit normal to $B$ at $w\in\partial B$. Considering a stochastic differential equation with normal reflection field on the boundary of $B$, the results of [8] (see also [10]) show that there exists a system $\Phi=(\Lambda,\mathcal{G},\{\mathcal{G}_t\},\mathbb{Q},Z)$ and continuous, adapted processes $W$ and $\ell$ such that, a.s., $W_t\in B$ for all $t\ge 0$, $\ell$ is continuous and locally of bounded variation,

$$W_t = w_0 + \int_0^t \vartheta(W_s)\,ds + \int_0^t \sigma(W_s)\,dZ_s + \ell_t, \qquad t\ge 0, \tag{2.9}$$

$$\ell_t = \int_0^t 1_{\{W_s\in\partial B\}} n(W_s)\,d|\ell|_s, \qquad t\ge 0. \tag{2.10}$$

Define a continuous, adapted process $U$ such that, a.s.,

$$U_t \doteq \int_0^t 1_{\{W_s\in\partial B\}}\varpi(n(W_s))\,d|\ell|_s.$$

It follows from Lemma 2.2 that (2.1) is satisfied. Also, Definition 2.1(iii) holds clearly, and $U$ satisfies part (ii) of Definition 2.1 since $\mathcal{U}$ is convex, and so $U\in\mathcal{A}(w_0)$. Let

$$\phi(w) = -(2\varepsilon)^{-1}|w-w_0|^2.$$

Note that $\phi$ is bounded on $B$ and satisfies $n\cdot\nabla\phi = 1$ on $\partial B$. Itô's formula applied to (2.9) shows that, a.s.,

$$\phi(W_t) = \phi(w_0) - \int_0^t \mathcal{L}\phi(W_s)\,ds + M_t + \int_0^t \nabla\phi(W_s)\cdot d\ell_s, \tag{2.11}$$

for $\mathcal{L}$ as in (2.8) and a continuous martingale

$$M = \left\{\int_0^t \nabla\phi(W_s)\cdot\sigma(W_s)\,dZ_s, t\ge 0\right\}$$

that starts from 0. Using (2.10) and the fact that $n \cdot \nabla\phi = 1$ on $\partial B$, the last term of (2.11) is equal to $|\ell|_t$, and therefore, $\mathbb{E}[|\ell|_t] \le C(1+t)$ for an appropriate constant $C$ not depending on $t$ since $\phi$ and $\mathcal{L}\phi$ are bounded on $B$. As a result, a similar bound holds for $\mathbb{E}[|U_t|]$, and it follows that $V(w_0) < \infty$ for $w_0 \in \mathcal{W}^o$.

We now consider Lipschitz continuity on the interior and finiteness on the boundary. The first property is essentially a consequence of (1.2) and Lemma 2.2 since they imply that the state process can be moved from any $w_1 \in \mathcal{W}$ to any other $w_2 \in \mathcal{W}$, instantaneously, by exercising a control and paying a cost that is bounded by a constant times $|w_1 - w_2|$. More precisely, fix $w_1 \in \mathcal{W}^o$ and $w_2 \in \mathcal{W}$. Given $\varepsilon \in (0,\infty)$, let $U \in \mathcal{A}(w_1)$ be such that $V(w_1) \le J(w_1, U) \le V(w_1) + \varepsilon$. Note that $U(t) \in \mathcal{U}$ for all $t \ge 0$, a.s. Let $\widetilde{u} \doteq \varpi(w_1 - w_2)$ and define $\widetilde{U} \doteq \widetilde{u} + U$. Clearly, $\widetilde{U} \in \mathcal{A}(w_2)$ (we use the fact that $\mathcal{U}$ is convex and a cone for this). Also,

$$V(w_2) - V(w_1) \le J(w_2, \widetilde{U}) - J(w_1, U) + \varepsilon = \kappa \cdot \widetilde{u} + \varepsilon \le c_\varpi |\kappa||w_1 - w_2| + \varepsilon.$$

Letting $\varepsilon \to 0$, it follows that the Lipschitz property holds on $\mathcal{W}^o$ and that $V$ is finite on the boundary. The argument above can now be repeated to deduce the Lipschitz property on all of $\mathcal{W}$. $\square$

Recall the notation $\mathcal{H}(q)$ from (1.4) and that $\mathcal{H}(q) \in (-\infty, +\infty)$ for every $q \in \mathbb{R}^d$. The following partial differential equation of Hamilton–Jacobi–Bellman type will be associated with our control problem:

$$((\mathcal{L} + \alpha)\psi - g) \vee \mathcal{H}(D\psi) = 0, \tag{2.12}$$

with an additional state constraint boundary condition. The precise definition of a solution for this equation with boundary condition is as a constrained viscosity solution of (2.12) on $\mathcal{W}$, defined in the following way (the form of the boundary condition was introduced by Soner in the paper [9]; that paper also contains an explanation on how the boundary condition was derived for the state constraint problem).

Definition 2.2 (Constrained viscosity solution).

(i) For $S \subset \mathcal{W}$, $\psi$ is said to be a viscosity supersolution [resp., subsolution] of (2.12) on $S$ if $\psi$ is a continuous real valued function on $\overline{S}$ and for all $w \in S$ and all $\varphi \in C^2(S)$ for which $\psi - \varphi$ has a global minimum [maximum] on $S$ at $w$, one has

$$(\mathcal{L}\varphi(w) + \alpha\psi(w) - g(w)) \vee \mathcal{H}(D\varphi(w)) \ge 0 \ [\le 0].$$

(ii) A function $\psi : \mathcal{W} \to \mathbb{R}$ is said to be a constrained viscosity solution of (2.12) on $\mathcal{W}$ if it is a viscosity subsolution of (2.12) on $\mathcal{W}^o$ and a viscosity supersolution of (2.12) on $\mathcal{W}$.

The latter condition in (ii) corresponds to the boundary condition. We remark that the definition above is equivalent to one in which the terms "global minimum" and "global maximum" are replaced by "local minimum" and, respectively, "local maximum," as is easy to see; this will be used in the sequel several times. Our main result characterizes the value function as a constrained viscosity solution of (2.12) on $\mathcal{W}$.

THEOREM 2.1. (i) *Solvability. $V$ is a constrained viscosity solution of (2.12) on $\mathcal{W}$.*

(ii) *Uniqueness. Condition (1.7) is necessary and sufficient for $V$ to be the only constrained viscosity solution.*

PROOF. Part (i) of the theorem follows from Propositions 4.2 and 4.3 below, whereas part (ii) is established in Corollary 5.1. □

The following simple example illustrates the effect of condition (1.7) on uniqueness.

EXAMPLE 2.1. (a) Let $\mathcal{W} = [0,1]$, $\mathcal{U} = \mathbb{R}_+^2$, $G = [1,-1]$, $\vartheta = 0$, $\sigma = 0$. Note that (1.2) is satisfied in this example. An admissible control $U = (U_1, U_2)'$ is nondecreasing in both coordinates and constrains $W(t)$ to $\mathcal{W}$ for all $t \geq 0$, where

$$W = w_0 + GU = w_0 + U_1 - U_2.$$

Let $g = 0$, $\alpha = 1$, and $\kappa = [\kappa_0, -\kappa_0]'$, $\kappa_0 \in (0,\infty)$. Then

$$\kappa \cdot U = \kappa_0(U_1 - U_2) = \kappa_0(W - w_0). \tag{2.13}$$

It is then easy to see that $V(w) = -\kappa_0 w$. Equation (2.12) on $(0,1)$ is equivalent to $\psi \vee \mathcal{H}(D\psi) = 0$, where

$$\mathcal{H}(q) = \sup\{-qGu - \kappa(u) : u \in \mathcal{U}_1\}, \qquad q \in \mathbb{R},$$

$\kappa(u) = \kappa_0(u_1 - u_2)$ and $\mathcal{U}_1 = \{u \in \mathbb{R}_+^2 : |u_1 - u_2| = 1\}$. Thus,

$$\mathcal{H}(q) = \sup\{(u_1 - u_2)(-q - \kappa_0) : u \in \mathbb{R}_+^2, |u_1 - u_2| = 1\} = |q + \kappa_0|.$$

In particular, (1.5) is satisfied. We show that $\psi(w) = -\kappa_0 w - c$ is a constrained viscosity solution of (2.12) on $\mathcal{W}$ for every $c \geq 0$. On $(0,1)$ the function $\psi$ satisfies the equation classically, and as is well known, this is sufficient for it to satisfy the equation in the viscosity sense [3] and hence to be a viscosity subsolution and supersolution of (2.12) on $\mathcal{W}^o$. Let us demonstrate that the boundary condition holds. If $\varphi \in C^2(\mathcal{W})$ is such that $\psi - \varphi$ has a global minimum at $w = 0$, then $D\varphi(0) = -\kappa_0 - c_1$, where $c_1 \geq 0$. The condition to be verified, $(-c) \vee |-c_1| \geq 0$, clearly holds. So does the condition to be verified at $w=1$, that is, $(-\kappa_0 - c) \vee |c_1| \geq 0$. Thus, $\psi$ is a constrained viscosity solution of (2.12) on $\mathcal{W}$ for every $c \geq 0$. The equation thus has multiple solutions. Note that in this example (1.7) is not satisfied.

(b) Consider now the case where $\kappa = [2\kappa_0, -\kappa_0]'$, $\kappa_0 \in (0,\infty)$. Once again, $V(w) = -\kappa_0 w$. Also note that, for $q \in \mathbb{R}$,

$$\begin{aligned}\mathcal{H}(q) &= \sup\{-(u_1 - u_2)q - \kappa_0(2u_1 - u_2) : u \in \mathbb{R}^2_+, |u_1 - u_2| = 1\}\\ &= \sup\{(u_1 - u_2)(-q - \kappa_0) - u_1\kappa_0 : u \in \mathbb{R}^2_+, |u_1 - u_2| = 1\}\\ &= (q + \kappa_0) \vee (-q - 2\kappa_0).\end{aligned}$$

Hence, $\mathcal{H}(\cdot)$ assumes both positive and negative values, and so (1.7) holds. Also, (2.12) for this problem can be written as

$$\psi \vee \mathcal{H}(D\psi) = \psi \vee (D\psi + \kappa_0) \vee (-D\psi - 2\kappa_0) = 0. \tag{2.14}$$

From the uniqueness statement of Theorem 2.1 it follows that $V$ is the unique constrained viscosity solution of (2.14). Also, note that each of the functions $\psi(w) = -\kappa_0 w - c$, $c \geq 0$ satisfies the equation classically on $(0,1)$, while only for $c = 0$ the boundary condition is satisfied. Indeed, for these functions, the boundary condition at $w = 0$ and, respectively, $w = 1$ reads $(-c) \vee (-c_1) \vee (c_1 - \kappa_0) \geq 0$ and $(-\kappa_0 - c) \vee c_1 \vee (-\kappa_0 - c_1) \geq 0$ for every $c_1 \geq 0$, which holds if and only if $c = 0$.

Next, we study some relations between the no arbitrage condition (1.8) and the hypotheses of Theorem 2.1. A condition similar to (1.8) was used in [7] in proving the existence of admissible controls and finiteness of the value function for the control problem therein. We note that there are cases covered by Proposition 2.1 and Theorem 2.1 in which the no arbitrage condition does not hold. This is the case in Example 2.1(a) in which, as mentioned above, our standing assumptions (1.2) and (1.5) hold, while (1.8) and (1.7) fail.

We now show that the no arbitrage condition is, in fact, useful in verifying (1.7). We introduce the following additional condition:

$$\text{there exists a unit vector } u_1 \in \mathbb{R}^p \text{ such that } \inf_{u \in \mathcal{U}_1} u_1 \cdot u > 0. \tag{2.15}$$

By (1.2), $\inf\{|u| : u \in \mathcal{U}, |Gu| = 1\}$ is strictly positive and, thus, a sufficient condition for (2.15) to hold is that there is a vector $u_1 \in \mathbb{R}^p$ for which

$$\inf_{u \in \mathcal{U} : |u| = 1} u_1 \cdot u > 0.$$

Indeed,

$$\inf_{u \in \mathcal{U}_1} u_1 \cdot u = \inf_{u \in \mathcal{U}_1} u_1 \cdot \frac{u}{|u|}|u| \geq \delta \inf_{u \in \mathcal{U} \setminus \{0\}} u_1 \cdot \frac{u}{|u|} = \delta \inf_{u \in \mathcal{U} : |u| = 1} u_1 \cdot u,$$

where $\delta = \inf\{|u| : u \in \mathcal{U}, |Gu| = 1\}$.

THEOREM 2.2. *Let conditions (1.8) and (2.15) hold. Then (1.7) holds.*

The proof of the theorem is presented below Corollary (2.1).

A weaker version of the no arbitrage condition (1.8) is (1.9). The latter condition is in fact equivalent to (1.5) (cf. Lemma A.1). We thus have two sets of sufficient conditions for existence and uniqueness of solutions of the DPE: In terms of the Hamiltonian $\mathcal{H}$ [(1.5) for existence and (1.7) for uniqueness], and in terms of no arbitrage considerations [(1.9) for existence and (1.8) and (2.15) for uniqueness].

REMARK 2.2. (a) Note that in Example 2.1(a), (2.15) is satisfied and, thus, from Theorem 2.2, (1.8) fails. However, the "weak no arbitrage" condition (1.9) holds in this example since it is equivalent to (1.5).

(b) Condition (1.8) is in fact strictly stronger [under (2.15)] than (1.7) as the following example shows. Let $\mathcal{W} = [0,1]$, $\mathcal{U} = \mathbb{R}^3_+$, $G = [1,-1,0]$, $\vartheta = 0$, $\sigma = 0$, $\kappa = [2\kappa_0, -\kappa_0, 0]'$. It is easy to see that $\mathcal{H}$ is the same as in Example 2.1(b) and thus (1.7) holds. However, (1.8) clearly fails.

The following corollary is an immediate consequence of the theorem.

COROLLARY 2.1. *Under the hypotheses of Theorem 2.2, $V$ is the only constrained viscosity solution to (2.12) on $\mathcal{W}$.*

PROOF OF THEOREM 2.2. For $\widetilde{\kappa} \in \mathbb{R}^p$, let

$$\mathcal{E}(\widetilde{\kappa}) \doteq \{u \in \mathcal{U} : Gu = 0 \text{ and } \widetilde{\kappa} \cdot u \le 0\}.$$

Note that $\mathcal{E}(\widetilde{\kappa})$ is always nonempty as it contains $\{0\}$. From condition (1.8) we have that $\mathcal{E}(\kappa) = \{0\}$. Next we show that $\{\widetilde{\kappa} \in \mathbb{R}^p : \mathcal{E}(\widetilde{\kappa}) = \{0\}\}$ is an open set. This is equivalent to showing that $F \doteq \{\widetilde{\kappa} \in \mathbb{R}^p : \mathcal{E}(\kappa) \ne \{0\}\}$ is closed. To this end, let $\kappa_n$ be a sequence of vectors in $F$ which converges to $\kappa^*$. Let $u_n \ne 0$ be in $\mathcal{E}(\kappa_n)$ for each $n$. Noting that the set $\mathcal{E}(\widetilde{\kappa})$ is a cone, we can assume, without loss of generality, that $|u_n| = 1$ for all $n$. Also, by choosing a subsequence if needed, we can assume that $u_n$ converges to some unit vector $u^* \in \mathcal{U}$. Note that $Gu^* = 0$. Since $\kappa_n \cdot u_n \to \kappa^* \cdot u^*$, we have $\kappa^* \cdot u^* \le 0$. This shows that $\mathcal{E}(\kappa^*) \ni u^* \ne 0$ and, as a result, $\kappa^* \in F$. Hence, $F$ is closed, and we thus have the following:

$$\text{there is } \varepsilon_0 > 0 \text{ such that } \mathcal{E}(\widetilde{\kappa}) = \{0\} \text{ whenever } |\widetilde{\kappa} - \kappa| < \varepsilon_0.$$

For $\varepsilon \in (0, \varepsilon_0)$, let $\kappa^\varepsilon \doteq \kappa - \varepsilon u_1$, where $u_1$ is as in (2.15). By the above display, $\mathcal{E}(\kappa^\varepsilon) = \{0\}$. Define $V^\varepsilon(\cdot)$ by (2.7) with $J(w_0, U)$ defined by (2.6) with $\kappa$ replaced by $\kappa^\varepsilon$. Note that condition (1.8) holds with $\kappa$ replaced by $\kappa^\varepsilon$. Consequently, by Lemma A.1 and since (1.8) implies (1.9), (1.5) holds with $\kappa$ replaced by $\kappa^\varepsilon$. In particular, applying Proposition 2.1 with $\kappa$ and $V$ replaced by $\kappa^\varepsilon$ and $V^\varepsilon$, we see that $V^\varepsilon(w) \in (-\infty, \infty)$ for all $w \in \mathcal{W}$ and $V^\varepsilon$ is Lipschitz continuous. Next, let $\zeta \in C^2(\mathcal{W})$ be a nonnegative function on $\mathcal{W}$ with $\zeta = 0$ on $\partial\mathcal{W}$ and $\sup_{w\in\mathcal{W}} \zeta(w) > 0$. Such a function can be easily constructed by choosing a nonnegative (but not identically zero) function that is twice continuously differentiable with compact support in the interior of $\mathcal{W}$. Let $\bar{M} = \max_{w\in\mathcal{W}} V^\varepsilon(w)$. Define

$$a = \inf\left\{\beta \geq 0 : \inf_{w\in\mathcal{W}} (\bar{M} + 1 - \beta\zeta(w) - V^\varepsilon(w)) \leq 0\right\}.$$

Clearly, $a$ is finite, $\inf_{w\in\mathcal{W}}(\bar{M} + 1 - a\zeta(w) - V^\varepsilon(w)) = 0$, and there exists $w \in \mathcal{W}^o$ such that $\bar{M} + 1 - a\zeta(w) - V^\varepsilon(w) = 0$. Define $\varphi \doteq \bar{M} + 1 - a\zeta$. Then $V^\varepsilon(w) = \varphi(w)$ and $V^\varepsilon(w_0) \leq \varphi(w_0)$ for all $w_0 \in \mathcal{W}$. For all $\delta > 0$ small enough, one has $w + \delta Gu \in \mathcal{W}^o$ for all $u \in \mathcal{U}_1$. Using Lemma A.3(i) in the Appendix, with $V$ replaced by $V^\varepsilon$, we obtain, for all $\delta > 0$ sufficiently small,

$$\varphi(w + \delta Gu) - \varphi(w) \geq V^\varepsilon(w + \delta Gu) - V^\varepsilon(w) \geq -\delta\kappa^\varepsilon \cdot u, \qquad u \in \mathcal{U}_1.$$

Dividing by $\delta$ and taking $\delta \to 0$, we have

$$D\varphi(w) \cdot Gu + \kappa \cdot u \geq \varepsilon u_1 \cdot u \qquad \text{for all } u \in \mathcal{U}_1.$$

Taking infimum over $u \in \mathcal{U}_1$ and using (2.15), we obtain $\mathcal{H}(q_0) < 0$ with $q_0 = D\varphi(w)$. $\square$

**3. Generalized Brownian networks.** Recently, in [7], Harrison and Williams considered a control problem for a Brownian network model and proved that it can be reduced to an equivalent but simpler control problem of a lower state dimension. In this section we note that the reduced control problem of [7] is a singular control problem with state constraints of the form introduced in Section 2 and prove that the standing assumptions of [7] imply those of the current paper (except for cases in which the reduced control problem is degenerate, i.e., the state space is a singleton). The data of a *Generalized Brownian Network* of [7] consists of the following:

(a) Positive integers $m, n, p$ specifying the dimensions of the state space, control space and control constraint space, respectively, for the Brownian network control problem.

(b) A vector $\theta \in \mathbb{R}^m$ and a nondegenerate $m \times m$ covariance matrix $\Sigma$.

(c) An $m \times n$ matrix $R$ and a $p \times n$ matrix $K$, which specify the effects of controls on the state of the system and constraints on the controls, respectively.

(d) A compact, convex set $\mathcal{Z} \subset \mathbb{R}^m$ that has a nonempty interior, which specifies the state space of the network control problem.

The cost process for the network control problem is specified by a continuous function $h : \mathcal{Z} \to \mathbb{R}$ and a vector $v \in \mathbb{R}^n$. The key assumptions in [7] that are made on the network data and the cost are the following.

ASSUMPTION 3.1 (*Assumption* 2.1 *of* [7]). $\{Ry : Ky \geq 0, y \in \mathbb{R}^n\} = \mathbb{R}^m$.

ASSUMPTION 3.2 (*Assumption* 2.2 *of* [7]).

$$\{y \in \mathbb{R}^n : Ky \geq 0, Ry = 0 \text{ and } v \cdot y \leq 0\} = \{0\}.$$

Lemma 4.3 of [7] shows that there is an $m$-dimensional vector $\pi$ and a $p$-dimensional vector $\kappa$ such that

$$v' = \pi' R + \kappa' K. \tag{3.1}$$

Define $\mathcal{M} \doteq \{a \in \mathbb{R}^m : a'R = b'K \text{ for some } b \in \mathbb{R}^p\}$. This is called the workload space in [7]. Let $d$ be the dimension of $\mathcal{M}$. Here we assume that $d > 0$, as we only treat state spaces for the reduced problem that have a nonempty interior. (The case $d = 0$ is a degenerate case that can occur; however, in this case the reduced control problem dramatically simplifies to one where the cost effectively only varies with the control and not with the state.) Let $M$ be a $d \times m$ matrix whose rows are a maximal linearly independent set of vectors in $\mathcal{M}$. Denote by $\mathcal{K}$ the range space of $K$. Define

$$\mathcal{W} \doteq \{Mz : z \in \mathcal{Z}\}, \qquad \mathcal{U} \doteq \mathcal{K} \cap \mathbb{R}^p_+, \qquad \vartheta \doteq M\theta, \qquad \Gamma \doteq M\Sigma M'. \tag{3.2}$$

Since $M$ is of full row rank, it follows that $\mathcal{W}$, like $\mathcal{Z}$, is compact and convex and has nonempty interior. Set $\sigma = \Gamma^{1/2}$, a positive definite square root of $\Gamma$. Note that $\vartheta$ and $\sigma$ are state independent. From Lemma 4.2 of [7] we have that there is a $d \times p$ matrix $G$ such that $MR = GK$. In this setting of a generalized Brownian network, we have the following.

LEMMA 3.1. *Assumption* 3.2 *implies condition (1.8).*

PROOF. Fix $u \in \mathcal{U}$ such that $Gu = 0$ and $\kappa \cdot u \leq 0$. From Lemma 4.4 of [7] (taking $x = 0$ therein) there exists a $y \in \mathbb{R}^n$ such that $u = Ky$, $Ry = 0$ and $v \cdot y = \kappa \cdot u$. From Assumption 3.2 we now have that $y = 0$. Thus, $u = Ky$ is zero as well and the lemma follows. □

Next, define the function $g : \mathcal{W} \to \mathbb{R}$ as

$$g(w) \doteq \inf\{h(z) + \alpha\pi \cdot z : Mz = w, z \in \mathcal{Z}\}, \qquad w \in \mathcal{W}. \tag{3.3}$$

The following continuous selection requirement is the final main assumption of [7].

Assumption 3.3 (*Assumption* 6.1 *of* [7]). There is a continuous function $\psi : \mathcal{W} \to \mathcal{Z}$ such that, for each $w \in \mathcal{W}$, $h(\psi(w)) + \alpha\pi \cdot \psi(w) = g(w)$ and $M\psi(w) = w$.

With the data $k = d$, $d$, $p$, $\mathcal{W}$, $\mathcal{U}$, $G$, $\vartheta$, $\sigma$, $g$ and $\kappa$, we now refer to our setting of Section 2: Admissible controls as defined in Definition 2.1 and the value function $V$ as defined in (2.7) comprise the reduced control problem of [7]. A minor difference from the setting in [7] is that there the control problem is formulated in terms of the process $\chi_t \doteq w_0 + \vartheta t + \sigma Z_t$ rather than in terms of $Z_t$, however, since given $w_0$ there is a one to one correspondence between $\chi$ and $Z$, the two formulations are easily seen to be equivalent. As a corollary of Theorem 2.1, we have the following.

Theorem 3.1. *Let Assumptions* 3.1, 3.2 *and* 3.3 *hold. Then $V$ is the unique constrained viscosity solution of* (2.12) *on $\mathcal{W}$.*

Proof. We first show that (1.2) holds. Fix $\xi \in \mathbb{R}^d$. Since the rows of $M$ are linearly independent, we can find $x \in \mathbb{R}^m$ such that $Mx = \xi$. From Assumption 3.1 we can find a $y \in \mathbb{R}^n$ such that $Ky \geq 0$ and $Ry = x$. Define $u = Ky$. Clearly, $u \in \mathcal{U}$. The result follows on noting that $Gu = GKy = MRy = Mx = \xi$.

Next, condition (1.8) holds by Lemma 3.1. Hence, (1.9) holds by Lemma A.2. Using the equivalence of (iii) and (v) of Lemma A.1 of the Appendix, we see that (1.5) holds as well. Finally, let $u_1 \in \mathbb{R}^p$ be a vector of which all entries are $p^{-1/2}$. Since $\mathcal{U} \subset \mathbb{R}^p_+$, $\inf_{u \in \mathcal{U}:|u|=1} u_1 \cdot u > 0$. As a result [see the comment below (2.15)] we also have $\inf_{u \in \mathcal{U}_1} u_1 \cdot u > 0$, and condition (2.15) follows. The result is now established by Corollary 2.1. □

**4. Dynamic programming principle and solvability.** In this section we will prove (i) of Theorem 2.1. We begin by introducing the following canonical representation for the control problem which facilitates the use of regular conditional probabilities in our proofs. Let $\mathcal{E} \doteq \mathcal{D}([0,\infty), \mathcal{W} \times \mathbb{R}^p) \times \mathcal{C}([0,\infty), \mathbb{R}^d)$. The canonical coordinate maps on $\mathcal{E}$ will be denoted by $\pi_i$, $i = 1,2,3$. For example, for $t \in [0,\infty)$, $\pi_1(t) : \mathcal{E} \to \mathcal{W}$ is defined as $\pi_1(\omega)(t) \doteq \omega_1(t)$ for $\omega \equiv (\omega_1, \omega_2, \omega_3) \in \mathcal{E}$. The definitions of $\pi_2, \pi_3$ are similar. Denote the triplet $(\pi_1, \pi_2, \pi_3)$ by $\pi$. We will endow $\mathcal{E}$ with the Borel sigma field $\mathcal{B}(\mathcal{E})$. Denote by $\mathcal{P}_{w_0}$ the class of all probability measures $\mathbb{P}^*$ on $(\mathcal{E}, \mathcal{B}(\mathcal{E}))$ satisfying the following:

- Under $\mathbb{P}^*$, $\pi_3$ is a $d$-dimensional $\{\mathcal{F}_t\}$-standard Brownian motion, where $\mathcal{F}_t$ is the canonical sigma field $\sigma\{\pi(s) : 0 \leq s \leq t\}$.

- $\pi_2$ is locally of bounded variation and has increments in $\mathcal{U}$, $\mathbb{P}^*$-a.s.
- The following holds, $\mathbb{P}^*$-a.s:

$$\pi_1(t) = w_0 + \int_0^t \vartheta(\pi_1(s))\,ds + \int_0^t \sigma(\pi_1(s))\,d\pi_3(s) + G\pi_2(t), \qquad t \geq 0, \tag{4.1}$$

where the second integral is an Itô integral;
-

$$\mathbb{E}^*\left[\int_0^\infty e^{-\alpha s}|\kappa\cdot\pi_2(s)|\,ds\right] < \infty, \tag{4.2}$$

where $\mathbb{E}^*$ is the expectation operator corresponding to $\mathbb{P}^*$.

In view of the fact that the value function is finite everywhere (Proposition 2.1), we may restrict to controls for which $J(w_0,U)$ is finite. Since $g$ is bounded on $\mathcal{W}$, by Lemma 2.1(i) and the equivalent representation of the cost in (2.6), we can represent $V$ as follows:

$$V(w_0) = \inf_{\mathbb{P}^*\in\mathcal{P}_{w_0}} \mathbb{E}^*\left[\int_0^\infty e^{-\alpha s}(g(\pi_1(s)) + \alpha\kappa\cdot\pi_2(s))\,ds\right].$$

Restricting to $\mathbb{P}^*$ that satisfy (4.2) considerably simplifies arguments in the proofs and the inequality (4.2) will be implicitly used at several places in the rest of this section. The following dynamic programming principle is key to the proof of part (i) of Theorem 2.1.

Proposition 4.1. *Fix $w_0\in\mathcal{W}$, and denote $\mathcal{W}^\varepsilon = \mathcal{W}\setminus\overline{B_\varepsilon(w_0)}$. Fix $\varepsilon>0$ such that $\mathcal{W}^\varepsilon$ is nonempty. Let*

$$\tau \doteq \inf\{s\geq 0 : \pi_1(s)\notin \overline{B_\varepsilon(w_0)}\}, \qquad \tau_t \doteq \tau\wedge t, t\in[0,\infty). \tag{4.3}$$

*Then for $t\in[0,\infty)$,*

$$\mathbb{E}^*(e^{-\alpha\tau_t}|\kappa\cdot\pi_2(\tau_t)|) < \infty \qquad \text{for all } \mathbb{P}^*\in\mathcal{P}_{w_0} \tag{4.4}$$

*and*

$$\begin{aligned} V(w_0) = \inf_{\mathbb{P}^*\in\mathcal{P}_{w_0}} \mathbb{E}^*\Bigg[&\int_0^{\tau_t} e^{-\alpha s}[g(\pi_1(s)) + \alpha\kappa\cdot\pi_2(s)]\,ds \\ &+ e^{-\alpha\tau_t}[V(\pi_1(\tau_t)) + \kappa\cdot\pi_2(\tau_t)]\Bigg]. \end{aligned} \tag{4.5}$$

Proof. Fix a $\mathbb{P}^*\in\mathcal{P}_{w_0}$ and $t\in[0,\infty)$. Let $\mathcal{G}_s = \mathcal{F}_{s+} = \bigwedge_{u>s}\mathcal{F}_u$ for all $s\in[0,\infty)$. Noting that $\tau_t$ is a $\{\mathcal{G}_s\}$-stopping time, apply Lemma A.4 with $(\Omega,\mathcal{F},P)$ in the lemma replaced by $(\mathcal{E},\mathcal{B}(\mathcal{E}),\mathbb{P}^*)$, $X$ replaced by $(\pi_1(\tau_t+\cdot),\pi_2(\tau_t+\cdot)-\pi_2(\tau_t),\pi_3(\tau_t+\cdot)-\pi_3(\tau_t))$, $\mathcal{T}$ replaced by $\mathcal{E}$ and $\mathcal{G}$ replaced by $\mathcal{G}_{\tau_t}$. We denote the resulting regular conditional probability measure by $\mathbb{P}^*_{\tau_t}$ and the associated conditional expectation operator by $\mathbb{E}^*_{\tau_t}$. For $\mathbb{P}^*$-a.e. $\omega$, $\mathbb{P}^*_{\tau_t}(\omega,\cdot)$ is a probability measure on $(\mathcal{E},\mathcal{B}(\mathcal{E}))$. Denote

$$E=\{\omega\in\mathcal{E}:\mathbb{P}^*_{\tau_t}(\omega,\cdot)\in\mathcal{P}_{\pi_1(\tau_t(\omega))}\}.$$

Observe that, $\mathbb{P}^*$-a.s., for $s\geq 0$,

$$\begin{aligned}\pi_1(\tau_t+s)=\ &\pi_1(\tau_t)+\int_{\tau_t}^{\tau_t+s}\vartheta(\pi_1(\eta))\,d\eta\\&+\int_{\tau_t}^{\tau_t+s}\sigma(\pi_1(\eta))\,d\pi_3(\eta)+G[\pi_2(\tau_t+s)-\pi_2(\tau_t)].\end{aligned}$$

Using the independence, under $\mathbb{P}^*$, of $\pi_3(\tau_t+\cdot)-\pi_3(\tau_t)$ from $\mathcal{G}_{\tau_t}$, as follows from the fact that $\pi_3$ is an $\{\mathcal{F}_t\}$-standard Brownian motion under $\mathbb{P}^*$, we have that for $\mathbb{P}^*$- a.e. $\omega$, under $\mathbb{P}^*_{\tau_t}(\omega,\cdot)$, $\pi_3$ is a $d$-dimensional $\{\mathcal{F}_t\}$-standard Brownian motion. Also, $\mathbb{P}^*$-a.s.,

$$\begin{aligned}\mathbb{E}^*_{\tau_t}\left[\int_0^\infty e^{-\alpha s}|\kappa\cdot\pi_2(s)|\,ds\right]&=\mathbb{E}^*\left[\int_0^\infty e^{-\alpha s}|\kappa\cdot(\pi_2(s+\tau_t)-\pi_2(\tau_t))|\,ds\Big|\mathcal{G}_{\tau_t}\right]\\&\leq\mathbb{E}^*\left[\int_0^\infty e^{-\alpha s}|\kappa\cdot\pi_2(s+\tau_t)|\,ds\Big|\mathcal{G}_{\tau_t}\right]\\&\quad+\alpha^{-1}|\kappa\cdot\pi_2(\tau_t)|\\&\leq e^{\alpha\tau_t}\mathbb{E}^*\left[\int_{\tau_t}^\infty e^{-\alpha s}|\kappa\cdot\pi_2(s)|\,ds\Big|\mathcal{G}_{\tau_t}\right]\\&\quad+\alpha^{-1}|\kappa\cdot\pi_2(\tau_t)|\\&<\infty,\end{aligned}$$

where the last inequality follows $\mathbb{P}^*$-a.s. from the fact that every $\mathbb{P}^*\in\mathcal{P}_{w_0}$ satisfies (4.2). In a similar fashion one can verify that the second and third bullets above Proposition 4.1 hold with $\mathbb{P}^*$ there replaced by $\mathbb{P}^*_{\tau_t}(\omega,\cdot)$ and $w_0$ replaced by $\pi_1(\tau_t)(\omega)$ for $\mathbb{P}^*$-a.e. $\omega$. Combining these observations, we have that $\mathbb{P}^*(E)=1$, and $\mathbb{P}^*$-a.s., that

$$\begin{aligned}V(\pi_1(\tau_t))&\leq\mathbb{E}^*_{\tau_t}\left[\int_0^\infty e^{-\alpha s}(g(\pi_1(s))+\alpha\kappa\cdot\pi_2(s))\,ds\right]\\&=\mathbb{E}^*\left[\int_0^\infty e^{-\alpha s}[g(\pi_1(s+\tau_t))+\alpha\kappa\cdot(\pi_2(s+\tau_t)-\pi_2(\tau_t))]\,ds\Big|\mathcal{G}_{\tau_t}\right]\\&=e^{\alpha\tau_t}\mathbb{E}^*\left[\int_{\tau_t}^\infty e^{-\alpha s}[g(\pi_1(s))+\alpha\kappa\cdot\pi_2(s)]\,ds\,\Big|\mathcal{G}_{\tau_t}\right]-\kappa\cdot\pi_2(\tau_t).\end{aligned}\tag{4.6}$$

It follows that $\mathbb{P}^*$-a.s.,

$$\begin{aligned}e^{-\alpha\tau_t}\kappa\cdot\pi_2(\tau_t)\leq\ &\mathbb{E}^*_{\tau_t}\left[\int_{\tau_t}^\infty e^{-\alpha s}(g(\pi_1(s))+\alpha\kappa\cdot\pi_2(s))\,ds\right]\\&-e^{-\alpha\tau_t}V(\pi_1(\tau_t)).\end{aligned}\tag{4.7}$$

Using (4.2) and the boundedness of $V$ on $\mathcal{W}$, we now have that $e^{-\alpha\tau_t}\kappa\cdot\pi_2(\tau_t)$ is bounded above by a $\mathbb{P}^*$-integrable random variable and, thus,

$$\mathbb{E}^*(e^{-\alpha\tau_t}(\kappa\cdot\pi_2(\tau_t))^+)<\infty \qquad \text{for all } \mathbb{P}^*\in\mathcal{P}_{w_0}. \tag{4.8}$$

In particular, for fixed $\mathbb{P}^*\in\mathcal{P}_{w_0}$, the expression on the right-hand side of (4.5) is well defined (possibly taking value $-\infty$). An argument later in this proof [see (4.13)] will in fact show that (4.8) holds with $(\kappa\cdot\pi_2(\tau_t))^+$ replaced by $|\kappa\cdot\pi_2(\tau_t)|$.

From (4.7) it also follows that the quantity

$$\mathbb{E}^*\left[\int_0^\infty e^{-\alpha s}(g(\pi_1(s))+\alpha\kappa\cdot\pi_2(s))\,ds\right]$$

is bounded below by

$$\mathbb{E}^*\left[\int_0^{\tau_t} e^{-\alpha s}(g(\pi_1(s))+\alpha\kappa\cdot\pi_2(s))\,ds+e^{-\alpha\tau_t}[V(\pi_1(\tau_t))+\kappa\cdot\pi_2(\tau_t)]\right].$$

Taking the infimum over all $\mathbb{P}^*\in\mathcal{P}_{w_0}$, we see that the left-hand side of (4.5) is bounded below by the right-hand side of the same.

Next we establish the reverse inequality. Let $\delta\in(0,\infty)$ be arbitrary. It is easy to see that there is a countable set $\Lambda_\delta\subset\mathcal{W}$ and a measurable map $\lambda_\delta:\mathcal{W}\to\Lambda_\delta$ such that $\lambda_\delta(\mathcal{W}^\varepsilon)\subset\mathcal{W}^\varepsilon$ and $\sup_{w\in\mathcal{W}}|\lambda_\delta(w)-w|<\delta$. Define $\hat{\lambda}_\delta(w)\doteq\lambda_\delta(w)-w$, $w\in\mathcal{W}$. For each $w\in\Lambda_\delta$, let $\hat{\mathbb{P}}_w\in\mathcal{P}_w$ be such that

$$\hat{\mathbb{E}}_w\left[\int_0^\infty e^{-\alpha s}[g(\pi_1(s))+\alpha\kappa\cdot\pi_2(s)]\,ds\right]\le V(w)+\delta, \tag{4.9}$$

where $\hat{\mathbb{E}}_w$ is the expectation operator corresponding to the measure $\hat{\mathbb{P}}_w$. Fix $\mathbb{P}^*\in\mathcal{P}_{w_0}$. Let $\mathcal{E}_2\doteq\mathcal{E}\times\mathcal{E}$ and let $\theta\doteq(\hat{\theta},\widetilde{\theta})$ denote the canonical coordinate maps with $\hat{\theta}\equiv(\hat{\theta}_1,\hat{\theta}_2,\hat{\theta}_3)$ and $\widetilde{\theta}\equiv(\widetilde{\theta}_1,\widetilde{\theta}_2,\widetilde{\theta}_3)$. More precisely, denoting a typical element of $\mathcal{E}_2$ by $\omega=(\hat{\omega},\widetilde{\omega})$, we have for $s\in[0,\infty)$, $\hat{\theta}_i(s)(\omega)=\hat{\omega}_i(s)$, $\widetilde{\theta}_i(s)(\omega)=\widetilde{\omega}_i(s)$, $i=1,2,3$. Define $\hat{\tau}_t(\hat{\omega},\widetilde{\omega})\doteq\tau_t(\hat{\omega})$. Then $\hat{\tau}_t$ is a $\{\mathcal{G}_t^1\}$ stopping time, where $\mathcal{G}_t^1\doteq\mathcal{F}_{t+}^1$ and $\mathcal{F}_t^1\doteq\sigma\{\hat{\theta}(s):s\le t\}$ for $t\ge 0$. Consider the probability measure $\mathbb{Q}_{w_0}$ on $(\mathcal{E}_2,\mathcal{B}(\mathcal{E}_2))$, that for $A,B\in\mathcal{B}(\mathcal{E})$ satisfies

$$\mathbb{Q}_{w_0}(A\times B)\doteq\int_A\hat{\mathbb{P}}_{\lambda_\delta(\hat{\theta}_1(\hat{\tau}_t))}(B)\,d\mathbb{P}^*(\hat{\omega}).$$

[Since $\lambda_\delta(\hat{\theta}_1(\hat{\tau}_t))$ takes only countably many values (in $\Lambda_\delta$), this is a valid measure theoretic construction.] Define stochastic processes $Z,U$ and $W$ on $(\mathcal{E}_2,\mathcal{B}(\mathcal{E}_2),\mathbb{Q}_{w_0})$ such that, for $s\ge 0$,

$$W(s)\doteq\hat{\theta}_1(s)1_{[0,\hat{\tau}_t)}(s)+\widetilde{\theta}_1(s-\hat{\tau}_t)1_{[\hat{\tau}_t,\infty)}(s),$$

$$U(s)\doteq\hat{\theta}_2(s)1_{[0,\hat{\tau}_t)}(s)+(\hat{\theta}_2(\hat{\tau}_t)+\varpi(\hat{\lambda}_\delta(\hat{\theta}_1(\hat{\tau}_t)))+\widetilde{\theta}_2(s-\hat{\tau}_t))1_{[\hat{\tau}_t,\infty)}(s),$$

$$Z(s)\doteq\hat{\theta}_3(s)1_{[0,\hat{\tau}_t)}(s)+(\widetilde{\theta}_3(s-\hat{\tau}_t)+\hat{\theta}_3(\hat{\tau}_t))1_{[\hat{\tau}_t,\infty)}(s),$$

where $\varpi$ is as in Lemma 2.2.

Denote the measure induced by $(W,U,Z)$ on $(\mathcal{E},\mathcal{B}(\mathcal{E}))$ by $\tilde{\mathbb{P}}_{w_0}$. It is easily seen that with $\tilde{\mathbb{P}}_{w_0}$ in place of $\mathbb{P}^*$, the first three bullet points of this section hold. Thus, from Lemma 2.1(i),

$$\tilde{\mathbb{E}}_{w_0}\left[\int_0^\infty e^{-\alpha s}\alpha(\kappa\cdot\pi_2(s))^-\,ds\right]\le \alpha L<\infty, \tag{4.10}$$

where $\tilde{\mathbb{E}}_{w_0}$ is the expectation operator corresponding to $\tilde{\mathbb{P}}_{w_0}$ and we have denoted by $L$ the finite quantity in Lemma 2.1(i). Next note that

$$\begin{aligned}&\tilde{\mathbb{E}}_{w_0}\left[\int_{\tau_t}^\infty \alpha e^{-\alpha s}(\kappa\cdot\pi_2(s))^-\,ds\right]\\ &\quad=\mathbb{E}_{\mathbb{Q}_{w_0}}\left[\int_{\hat{\tau}_t}^\infty \alpha e^{-\alpha s}(\kappa\cdot[\hat{\theta}_2(\hat{\tau}_t)+\varpi(\hat{\lambda}_\delta(\hat{\theta}_1(\hat{\tau}_t)))+\widetilde{\theta}_2(s-\hat{\tau}_t)])^-\,ds\right].\end{aligned} \tag{4.11}$$

From the above equality and (4.10) we have that

$$\begin{aligned}&\mathbb{E}_{\mathbb{Q}_{w_0}}(e^{-\alpha\hat{\tau}_t}(\kappa\cdot\hat{\theta}_2(\hat{\tau}_t))^-)\\ &\quad\le \alpha L+c_\varpi|\kappa|\delta+\mathbb{E}_{\mathbb{Q}_{w_0}}\left[\int_{\hat{\tau}_t}^\infty \alpha e^{-\alpha s}(\kappa\cdot\widetilde{\theta}_2(s-\hat{\tau}_t))^-\,ds\right]\\ &\quad\le 2\alpha L+c_\varpi|\kappa|\delta<\infty,\end{aligned} \tag{4.12}$$

where the first inequality uses Lemma 2.2 and the second inequality uses once more (4.10) and Lemma 2.1(i). Thus, we have shown that

$$\mathbb{E}^*(e^{-\alpha\tau_t}(\kappa\cdot\pi_2(\tau_t))^-)<\infty\qquad\text{for all }\mathbb{P}^*\in\mathcal{P}_{w_0}. \tag{4.13}$$

Combining (4.8) and (4.13), we have (4.4) and, consequently,

$$\tilde{\mathbb{E}}_{w_0}\left[\int_0^\infty e^{-\alpha s}|\kappa\cdot\pi_2(s)|\,ds\right]<\infty.$$

In particular, recalling the properties of $\tilde{\mathbb{P}}_{w_0}$ stated above (4.10), we have that $\tilde{\mathbb{P}}_{w_0}\in\mathcal{P}_{w_0}$. Next,

$$\begin{aligned}&\tilde{\mathbb{E}}_{w_0}\left[\int_0^{\tau_t} e^{-\alpha s}[g(\pi_1(s))+\alpha\kappa\cdot\pi_2(s)]\,ds\right]\\ &\quad=\mathbb{E}^*\left[\int_0^{\tau_t} e^{-\alpha s}[g(\pi_1(s))+\alpha\kappa\cdot\pi_2(s)]\,ds\right].\end{aligned} \tag{4.14}$$

Also,

$$\begin{aligned}&\tilde{\mathbb{E}}_{w_0}\left[\int_{\tau_t}^\infty e^{-\alpha s}[g(\pi_1(s))+\alpha\kappa\cdot\pi_2(s)]\,ds\right]\\ &\quad=\mathbb{E}_{\mathbb{Q}_{w_0}}\left[\int_{\hat{\tau}_t}^\infty e^{-\alpha s}\{g(\widetilde{\theta}_1(s-\hat{\tau}_t))\right.\\ &\qquad\qquad\left.+\alpha\kappa\cdot[\hat{\theta}_2(\hat{\tau}_t)+\varpi(\hat{\lambda}_\delta(\hat{\theta}_1(\hat{\tau}_t)))+\widetilde{\theta}_2(s-\hat{\tau}_t)]\}\,ds\right]\end{aligned} \tag{4.15}$$

$$
\begin{aligned}
&= \mathbb{E}_{\mathbb{Q}_{w_0}}\left[e^{-\alpha\hat{\tau}_t}\int_0^\infty e^{-\alpha s}[g(\widetilde{\theta}_1(s)) + \alpha\kappa\cdot\widetilde{\theta}_2(s)]\,ds\right] \\
&\quad + \mathbb{E}_{\mathbb{Q}_{w_0}}[e^{-\alpha\hat{\tau}_t}\kappa\cdot[\hat{\theta}_2(\hat{\tau}_t) + \varpi(\hat{\lambda}_\delta(\hat{\theta}_1(\hat{\tau}_t)))]] \\
&\equiv T_1 + T_2.
\end{aligned}
$$

In splitting the expectation in the second equality above, we have used (4.4). As an immediate consequence of Lemma 2.2, we obtain

$$
T_2 \le \mathbb{E}^*[e^{-\alpha\tau_t}\kappa\cdot\pi_2(\tau_t)] + c_\varpi|\kappa|\delta. \tag{4.16}
$$

Note that by the definition of $\mathbb{Q}_{w_0}$, $T_1$ can be rewritten as

$$
\mathbb{E}^*\left[e^{-\alpha\tau_t}\hat{\mathbb{E}}_{\lambda_\delta(\pi_1(\tau_t))}\left\{\int_0^\infty e^{-\alpha s}[g(\pi_1(s)) + \alpha\kappa\cdot\pi_2(s)]\,ds\right\}\right].
$$

From (4.9) and Proposition 2.1 we can bound the last expression from above by

$$
\mathbb{E}^*[e^{-\alpha\tau_t}V(\lambda_\delta(\pi_1(\tau_t)))] + \delta \le \mathbb{E}^*[e^{-\alpha\tau_t}V(\pi_1(\tau_t))] + (c_{\mathrm{lip}} + 1)\delta. \tag{4.17}
$$

Combining the above inequality with (4.15), (4.16) and (4.14), we obtain

$$
\begin{aligned}
V(w_0) &\le \tilde{\mathbb{E}}_{w_0}\left[\int_0^\infty e^{-\alpha s}[g(\pi_1(s)) + \alpha\kappa\cdot\pi_2(s)]\,ds\right] \\
&\le \mathbb{E}^*\left[\int_0^{\tau_t} e^{-\alpha s}[g(\pi_1(s)) + \alpha\kappa\cdot\pi_2(s)]\,ds\right. \\
&\qquad\qquad \left. + e^{-\alpha\tau_t}[V(\pi_1(\tau_t)) + \kappa\cdot\pi_2(\tau_t)]\right] \\
&\qquad\qquad + (1 + c_{\mathrm{lip}} + c_\varpi|\kappa|)\delta.
\end{aligned}
$$

Letting $\delta\to 0$ and taking infimum over all $\mathbb{P}^*\in\mathcal{P}_{w_0}$, we obtain the desired reverse inequality. □

PROPOSITION 4.2. *$V$ is a viscosity supersolution of (2.12) on $\mathcal{W}$.*

PROOF. Fix $w_0\in\mathcal{W}$ and let $\varphi\in C^2(\mathcal{W})$ be such that $V-\varphi$ has a global minimum at $w_0$. We can assume without loss of generality that $V(w_0) - \varphi(w_0) = 0$. We need to show that either

$$
\mathcal{L}\varphi(w_0) + \alpha\varphi(w_0) - g(w_0) \ge 0 \tag{4.18}
$$

or

$$
\inf\{Gu\cdot D\varphi(w_0) + \kappa\cdot u : u\in\mathcal{U}_1\} \le 0. \tag{4.19}
$$

Arguing by contradiction, assume that neither of the above assertions is true. Then one can find $\gamma>0$ and $\varepsilon>0$ such that, for all $\bar{w}\in\overline{B_{2\varepsilon}(w_0)}\cap\mathcal{W}$,

$$
\mathcal{L}\varphi(\bar{w}) + \alpha\varphi(\bar{w}) - g(\bar{w}) \le -\gamma \tag{4.20}
$$

and $Gu \cdot D\varphi(\bar{w}) + \kappa \cdot u \geq \gamma$ for all $u \in \mathcal{U}_1$. Note that the latter implies that

$$Gu \cdot D\varphi(\bar{w}) + \kappa \cdot u \geq \gamma |Gu| \qquad \text{for all } u \in \mathcal{U}. \tag{4.21}$$

Indeed, if $|Gu| > 0$, this is immediate. If $Gu = 0$, let $\bar{u}$ be a vector in $\mathcal{U}$ such that $|G\bar{u}| = 1$ [$\bar{u}$ exists by (1.2)]. Then for $r > 0$, $u^r \doteq u + r\bar{u} \in \mathcal{U}$ and $|Gu^r| > 0$, and so (4.21) holds for $u^r$ in place of $u$, and sending $r \to 0$, it follows that (4.21) holds for $u$ as well.

Let $t > 0$ and fix $\mathbb{P}^* \in \mathcal{P}_{w_0}$. Let $\varepsilon$, $\tau$, $\tau_t$ be as in Proposition 4.1. An application of Itô's formula gives

$$\begin{aligned}
\varphi(w_0) = {} & \mathbb{E}^*[e^{-\alpha\tau_t}\varphi(\pi_1(\tau_t))] \\
& + \mathbb{E}^*\left[\int_0^{\tau_t} e^{-\alpha s}[\mathcal{L}\varphi(\pi_1(s)) + \alpha\varphi(\pi_1(s))]\, ds\right] \\
& - \mathbb{E}^*\Bigg[\int_{[0,\tau_t]} e^{-\alpha s} D\varphi(\pi_1(s)) \cdot d\eta_s^c \\
& \qquad + \sum_{0 \leq s \leq \tau_t} e^{-\alpha s}[\varphi(\pi_1(s)) - \varphi(\pi_1(s-))]\Bigg],
\end{aligned} \tag{4.22}$$

where $\eta^c(s) \doteq G\varrho^c(s)$, and $\varrho^c(s) \doteq (\pi_2(s) - \sum_{0 \leq r \leq s} \Delta\pi_2(r))$ for $s \geq 0$, and $\pi_1(0-) = \omega_0$, $\Delta\pi_2(0) = \pi_2(0)$. Let

$$\mathcal{P}_{w_0,\varepsilon} \doteq \{\mathbb{P}^* \in \mathcal{P}_{w_0} : \mathbb{P}^*\{\tau < \infty, \pi_1(\tau) \notin \overline{B_{2\varepsilon}(w_0)}\} = 0\}.$$

From (4.21) we obtain, for $\mathbb{P}^* \in \mathcal{P}_{w_0,\varepsilon}$ and $s \in [0, \tau_t]$,

$$\begin{aligned}
\varphi(\pi_1(s)) - \varphi(\pi_1(s-)) &= \int_0^1 D\varphi(\pi_1(s-) + r\Delta(G\pi_2(s))) \cdot G\Delta\pi_2(s)\, dr \\
&\geq \gamma|\Delta(G\pi_2)(s)| - \kappa \cdot \Delta\pi_2(s), \qquad \mathbb{P}^*\text{-a.s.}
\end{aligned} \tag{4.23}$$

Note that in (4.23), $\pi_1(s-) + r\Delta(G\pi_2(s)) \in \mathcal{W}$ for all $r \in [0,1]$ since $\mathcal{W}$ is convex. Since $\pi_2$ has increments in $\mathcal{U}$, it is elementary to check that so does its continuous part $\varrho^c$. As another consequence of (4.21), we have, $\mathbb{P}^*$-a.s. (for $\mathbb{P}^* \in \mathcal{P}_{w_0,\varepsilon}$),

$$\int_{[0,\tau_t]} e^{-\alpha s} D\varphi(\pi_1(s)) \cdot d\eta_s^c \geq \gamma e^{-\alpha t}|\eta^c|_{\tau_t} - \int_{[0,\tau_t]} e^{-\alpha s}\, d(\kappa \cdot \varrho_s^c). \tag{4.24}$$

Combining the above two inequalities, we obtain

$$\begin{aligned}
\int_{[0,\tau_t]} & e^{-\alpha s} D\varphi(\pi_1(s)) \cdot d\eta_s^c + \sum_{0 \leq s \leq \tau_t} e^{-\alpha s}(\varphi(\pi_1(s)) - \varphi(\pi_1(s-))) \\
& \geq \gamma e^{-\alpha t}|G\pi_2|_{\tau_t} - \int_{[0,\tau_t]} e^{-\alpha s}\, d(\kappa \cdot \pi_2(s)) \\
& = \gamma e^{-\alpha t}|G\pi_2|_{\tau_t} - \alpha \int_0^{\tau_t} e^{-\alpha s}\kappa \cdot \pi_2(s)\, ds - e^{-\alpha\tau_t}\kappa \cdot \pi_2(t),
\end{aligned}$$

where the last equality follows from the integration by parts formula (2.2).

Using the above inequality and (4.20) in (4.22), we obtain

$$\begin{aligned}\varphi(w_0) \leq\ & \mathbb{E}^*[e^{-\alpha\tau_t}\varphi(\pi_1(\tau_t))] \\ & + \mathbb{E}^*\left[\int_0^{\tau_t} e^{-\alpha s}[g(\pi_1(s)) + \alpha\kappa\cdot\pi_2(s)]\,ds + e^{-\alpha\tau_t}\kappa\cdot\pi_2(t)\right] \\ & - \gamma e^{-\alpha t}\mathbb{E}^*[\tau_t + |G\pi_2|_{\tau_t}].\end{aligned}$$

Once again, in splitting the above expectation, we have used (4.2) and (4.4). Taking the infimum over all $\mathbb{P}^* \in \mathcal{P}_{w_0,\varepsilon}$, in the above inequality, we have

$$\varphi(w_0) \leq \inf_{\mathbb{P}^*\in\mathcal{P}_{w_0,\varepsilon}} \hat{J}(\mathbb{P}^*) - \gamma e^{-\alpha t}\beta(t), \tag{4.25}$$

where

$$\begin{aligned}\hat{J}(\mathbb{P}^*) \doteq\ & \mathbb{E}^*[e^{-\alpha\tau_t}V(\pi_1(\tau_t))] \\ & + \mathbb{E}^*\left[\int_0^{\tau_t} e^{-\alpha s}g(\pi_1(s))\,ds + \int_{[0,\tau_t]} e^{-\alpha s}\,d(\kappa\cdot\pi_2(s))\right],\end{aligned} \tag{4.26}$$

$\beta(t) \doteq \inf_{\mathbb{P}^*\in\mathcal{P}_{w_0,\varepsilon}} \mathbb{E}^*[\tau_t + |G\pi_2|_{\tau_t}]$, and we have used the fact that $\varphi \leq V$ in $\mathcal{W}$ plus the integration by parts formula (2.2). Using Lemma A.3 in the Appendix, the infimum on the right-hand side of (4.25) can be replaced by the infimum over all of $\mathcal{P}_{w_0}$. Thus, applying Proposition 4.1, we obtain $\varphi(w_0) \leq V(w_0) - \gamma e^{-\alpha t}\beta(t)$. Finally, to arrive at a contradiction, we show that

$$\text{there exists } t > 0 \text{ such that } \beta(t) > 0. \tag{4.27}$$

Fix $\mathbb{P}^* \in \mathcal{P}_{w_0,\varepsilon}$. Then $|\pi_1(\tau_t) - w_0| \geq \varepsilon$ on $\{\tau \leq t\}$, $\mathbb{P}^*$ a.s. Thus, denoting $r_t = \int_0^t \sigma(\pi_1(s))\,d\pi_3(s)$ and $|\vartheta|_{\mathcal{W}} = \max_{x\in\mathcal{W}}|\vartheta(x)|$,

$$\mathbb{E}^*[|G\pi_2|_{\tau_t}1_{\{\tau\leq t\}}] \geq \mathbb{E}^*[|G\pi_2(\tau_t)|1_{\{\tau\leq t\}}] \geq \mathbb{E}^*[(\varepsilon - t|\vartheta|_{\mathcal{W}} - |r|_t^*)1_{\{\tau\leq t\}}].$$

Hence, for $t \in (0, (3|\vartheta|_{\mathcal{W}})^{-1}\varepsilon)$,

$$\begin{aligned}\mathbb{E}^*[\tau_t + |G\pi_2|_{\tau_t}] &\geq t\mathbb{P}^*(\tau > t) + (\varepsilon/3)\mathbb{P}(\tau \leq t,\ |r|_t^* < \varepsilon/3) \\ &\geq [t \wedge (\varepsilon/3)]\mathbb{P}^*(|r|_t^* < \varepsilon/3).\end{aligned}$$

Clearly, for all $t > 0$ small enough, $\mathbb{P}^*(|r|_t^* < \varepsilon/3) > 0$. This proves (4.27) and hence the result. □

Next, we prove the subsolution property of the value function.

Proposition 4.3. *$V$ is a viscosity subsolution of (2.12) on $\mathcal{W}^o$.*

PROOF. Fix $w_0 \in \mathcal{W}^o$ and let $\varphi \in C^2(\mathcal{W}^o)$ be such that $V - \varphi$ has a global maximum on $\mathcal{W}^o$ at $w_0$. Once more, we can assume without loss of generality that $\varphi(w_0) = V(w_0)$. Thus, $V \le \varphi$ on $\mathcal{W}^o$. We need to show that $\zeta(w) \doteq \alpha\varphi(w) + \mathcal{L}\varphi(w) - g(w)$, $w \in \mathcal{W}^o$ satisfies

$$\zeta(w_0) \le 0 \tag{4.28}$$

and

$$Gu \cdot D\varphi(w_0) + \kappa \cdot u \ge 0, \qquad u \in \mathcal{U}_1. \tag{4.29}$$

For all $\delta > 0$ small enough, one has $w_0 + \delta G u \in \mathcal{W}^o$ for all $u \in \mathcal{U}_1$. By Lemma A.3,

$$\varphi(w_0 + \delta Gu) - \varphi(w_0) \ge V(w_0 + \delta Gu) - V(w_0) \ge -\delta\kappa \cdot u, \qquad u \in \mathcal{U}_1.$$

Dividing by $\delta$ and taking $\delta \to 0$ proves (4.29).

Now consider (4.28). Let $\varepsilon > 0$ be such that $\overline{B}_\varepsilon(w_0) \subset \mathcal{W}^o$. Let

$$\tau^\varepsilon \doteq 1 \wedge \inf\{t \ge 0 : \pi_1(t) \notin \overline{B}_\varepsilon(w_0)\}.$$

Now let $\mathbb{P}^* \in \mathcal{P}_{w_0}$ be such that $\mathbb{P}^*\{\pi_2(t) = 0, t \in [0, \tau^\varepsilon]\} = 1$. Such $\mathbb{P}^*$ exists by standard results on existence of solutions to the SDE (4.1) with $\pi_2 = 0$ and by the controllability condition (1.2) (for the behavior after $\tau^\varepsilon$). An application of Itô's formula gives

$$\begin{aligned} V(w_0) = \varphi(w_0) = {} & \mathbb{E}^*[e^{-\alpha\tau^\varepsilon}\varphi(\pi_1(\tau^\varepsilon))] \\ & + \mathbb{E}^*\left[\int_0^{\tau^\varepsilon} e^{-\alpha s}[\mathcal{L}\varphi(\pi_1(s)) + \alpha\varphi(\pi_1(s))]\,ds\right]. \end{aligned} \tag{4.30}$$

Using Proposition 4.1, $V \le \varphi$, and the above, we obtain

$$\begin{aligned} V(w_0) & \le \mathbb{E}^*\left[\int_0^{\tau^\varepsilon} e^{-\alpha s} g(\pi_1(s))\,ds + e^{-\alpha\tau^\varepsilon}\varphi(\pi_1(\tau^\varepsilon))\right] \\ & = \varphi(w_0) + \mathbb{E}^*\left[\int_0^{\tau^\varepsilon} e^{-\alpha s}[g(\pi_1(s)) - \mathcal{L}\varphi(\pi_1(s)) - \alpha\varphi(\pi_1(s))]\,ds\right]. \end{aligned} \tag{4.31}$$

Recalling that $V(w_0) = \varphi(w_0)$, we have $\mathbb{E}^*[\int_0^{\tau^\varepsilon} e^{-\alpha s}\zeta(\pi_1(s))\,ds] \le 0$. Hence,

$$\begin{aligned} \zeta(w_0)\mathbb{E}^*\left[\int_0^{\tau^\varepsilon} e^{-\alpha s}\,ds\right] & \le \mathbb{E}^*\left[\int_0^{\tau^\varepsilon} e^{-\alpha s}(\zeta(w_0) - \zeta(\pi_1(s)))\,ds\right] \\ & \le \alpha(w_0, \varepsilon)\mathbb{E}^*\left[\int_0^{\tau^\varepsilon} e^{-\alpha s}\,ds\right], \end{aligned}$$

where

$$\alpha(w_0, \varepsilon) = \max_{w \in \overline{B}_\varepsilon(w_0)} |\zeta(w) - \zeta(w_0)|.$$

Since $\tau^\varepsilon > 0$ $\mathbb{P}^*$-a.s., it follows that $\zeta(w_0) \le \alpha(w_0, \varepsilon)$. Taking $\varepsilon \to 0$, we obtain (4.28) by the continuity of $\zeta$ on $\mathcal{W}^o$. $\square$

**5. Uniqueness.** In this section we will prove part (ii) of Theorem 2.1. It will be a special case of a result that involves more general convex $\mathcal{H}$ that we now formulate. The state space $\mathcal{W}$ is assumed to be a closed, bounded subset of $\mathbb{R}^d$ satisfying the condition

$$\begin{aligned} &\text{for every } \xi\in\mathcal{W}, \text{ there exist } \eta=\eta(\xi)\in\mathbb{R}^d \text{ and } a=a(\xi)>0 \text{ such that}\\ &\quad B_{ta}(w+t\eta)\subset\mathcal{W}^o, \text{ for all } w\in\mathcal{W}\cap B_a(\xi) \text{ and all } t\in(0,1]. \end{aligned} \tag{5.1}$$

As shown in Lemma 6.1 of [1], any compact convex set with nonempty interior meets this condition. It is also satisfied by the closure of any bounded Lipschitz domain (in the sense of [2], Chapter III). Indeed, let $\mathcal{W}$ be the closure of a Lipschitz domain. For $\xi\in\mathcal{W}^o$, (5.1) obviously holds. Next, given $\xi\in\partial\mathcal{W}$, there is $a_1>0$, a Lipschitz function $f:\mathbb{R}^{d-1}\to\mathbb{R}$, and a coordinate system $CS$ such that $\mathcal{W}^o\cap B_{a_1}(\xi)=\{y\in B_{a_1}(\xi): y_1>f(y_2,\dots,y_d)\}$, where, for $y\in\mathbb{R}^d$, $(y_1,\dots,y_d)$ represents $y$ in the coordinate system $CS$. It is not hard to check that for $a>0$ small enough (depending only on $a_1$ and the Lipschitz constant of $f$) and $\eta=ae_1$ (in $CS$), one has $y_1>f(y_2,\dots,y_d)$ whenever $y=w+t\eta+taz$, $t\in(0,1]$, $z\in\mathbb{R}^d$, $|z|<1$, $w\in\mathcal{W}\cap B_a(\xi)$. Thus, (5.1) holds.

Next, we allow $\alpha$ and $\mathcal{H}$ to depend on $x$, and assume that for some $C\in[1,\infty)$ and a function $\omega:\mathbb{R}_+\to\mathbb{R}_+$ with vanishing right limit at 0 [i.e., $\omega(0+)=0$], one has

$$|\alpha(x)-\alpha(y)|+|\vartheta(x)-\vartheta(y)|+|\sigma(x)-\sigma(y)|\le C|x-y|,\qquad x,y\in\mathcal{W}, \tag{5.2}$$

$$\alpha(x)\ge\alpha_0>0,\qquad x\in\mathcal{W}, \tag{5.3}$$

$$\begin{aligned} &|g(x)-g(y)|+|\mathcal{H}(x,q)-\mathcal{H}(y,q)|\le\omega(|x-y|),\\ &\qquad\qquad x,y\in\mathcal{W},\ q\in\mathbb{R}^d. \end{aligned} \tag{5.4}$$

$$|\mathcal{H}(x,q_1)-\mathcal{H}(x,q_2)|\le C|q_1-q_2|,\qquad x\in\mathcal{W},\ q_1,q_2\in\mathbb{R}^d. \tag{5.5}$$

Recall that $\Gamma=\sigma\sigma'$. The constant $C\ge1$ will be assumed to be large enough so that

$$|\alpha(x)|+|\vartheta(x)|+|\Gamma(x)|+|\sigma(x)|\le C,\qquad x\in\mathcal{W}. \tag{5.6}$$

It is also assumed that $q\mapsto\mathcal{H}(x,q)$ is convex for every $x\in\mathcal{W}$. For $n\in\mathbb{N}$, denote by $\mathcal{S}(n)$ the space of symmetric $n\times n$ matrices, and write "$\le$" for the usual order on $\mathcal{S}(n)$ [for $A,B\in\mathcal{S}(n)$, $A\le B$ if and only if $B-A$ is nonnegative definite]. For $x\in\mathcal{W}$, $r\in\mathbb{R}$, $q\in\mathbb{R}^d$ and $A\in\mathcal{S}(d)$ denote

$$F(x,r,q,A)=\alpha(x)r-\vartheta(x)\cdot q-\tfrac12\operatorname{trace}(\Gamma(x)A)-g(x).$$

The comparison result below regards solutions to the equation

$$F(x,\psi,D\psi,D^2\psi)\vee\mathcal{H}(x,D\psi)=0, \tag{5.7}$$

defined with state constraint boundary condition (analogous to Definition 2.2). Since the left-hand side of (5.7) is not strictly increasing in the $\psi$ variable (because the term $\mathcal{H}$ does not depend on this variable), the argument for comparison will rely on the existence of appropriate strict subsolutions (cf. Section 5C of [3]). We therefore assume the following:

$$\begin{gathered}\text{there exists a constant } \gamma>0 \text{ and a function } \psi_0\in C^2(\mathcal{W}) \text{ such that}\\ \mathcal{H}(x,D\psi_0(x))\le-\gamma, \qquad x\in\mathcal{W}^o.\end{gathered} \tag{5.8}$$

THEOREM 5.1. *Let conditions (5.1)–(5.5), (5.8) hold, and let $\check v$ be a subsolution of (5.7) on $\mathcal{W}^o$ and let $v$ be a supersolution of (5.7) on $\mathcal{W}$. Then*

$$\check v\le v, \qquad \text{on } \mathcal{W}. \tag{5.9}$$

In the special case where $\mathcal{H}(x,q)$ is independent of $x$, condition (5.8) is easily seen to be equivalent to (1.7). Also, for this case, if $\inf_{q\in\mathbb{R}^d}\mathcal{H}(q)\ge0$, then the problem degenerates in the sense that every continuous function is a viscosity supersolution, and moreover, from any viscosity subsolution $\psi$, we can produce subsolutions $\psi-c$, $c>0$. These comments are summarized in the following.

COROLLARY 5.1. *Assume $\mathcal{H}(x,q)=\mathcal{H}(q)$, $x\in\mathcal{W}$, $q\in\mathbb{R}^d$, and let conditions (5.1)–(5.5) hold. Let there exist a constrained viscosity solution to (5.7). Then $\inf_{q\in\mathbb{R}^d}\mathcal{H}(q)<0$ is necessary and sufficient for uniqueness of such solutions.*

PROOF OF THEOREM 5.1. We introduce some notation specific to the proof. Let $S$ be a relatively open subset of $\mathcal{W}$. For $x\in S\subset\mathcal{W}$ and a real valued continuous function $\psi$ on $\mathcal{W}$, let the corresponding second order superjet and subjet be defined as follows [here we follow the terminology and notation of [3]; these objects are not to be confused with the cost functional $J$ of (1.3)]:

$$\begin{aligned}J_S^{2,+}\psi(x)&\doteq\{(D\varphi(x),D^2\varphi(x)):\\ &\qquad \varphi\in C^2(S) \text{ and } \psi-\varphi \text{ has a local maximum at } x\},\\ J_S^{2,-}\psi(x)&\doteq\{(D\varphi(x),D^2\varphi(x)):\\ &\qquad \varphi\in C^2(S) \text{ and } \psi-\varphi \text{ has a local minimum at } x\}.\end{aligned}$$

Define the closures of the above sets in the following way. For $x \in S$,

$$\begin{aligned}\overline{J}^{2,+}_S\psi(x) \doteq \{(q,M) \in \mathbb{R}^d \times \mathcal{S}(d):&\\ &\text{there exists a sequence } (x_n,q_n,M_n) \in \mathcal{W} \times \mathbb{R}^d \times \mathcal{S}(d) \text{ s.t.}\\ &(q_n,M_n) \in J^{2,+}_S\psi(x_n), \text{ and}\\ &\qquad (x_n,\psi(x_n),q_n,M_n) \to (x,\psi(x),q,M) \text{ as } n\to\infty\}.\end{aligned}$$

Define $\overline{J}^{2,-}_S\psi(x)$ analogously. For short, write $J^{2,+}$ for $J^{2,+}_{\mathcal{W}^o}$ and similarly define $J^{2,-}$, $\overline{J}^{2,+}$ and $\overline{J}^{2,-}$.

Let $\gamma$ and $\psi_0$ be as in (5.8) and set $\psi_1 = \psi_0 - c_1$, where the constant $c_1$ is large enough to ensure that

$$\begin{aligned}&F(x,\psi_1(x),D\psi_1(x),D^2\psi_1(x)) \le -1\\ &\text{and} \quad \psi_1(x) - \check{v}(x) \le 0 \qquad \text{for all } x \in \mathcal{W}.\end{aligned} \tag{5.10}$$

By (5.8),

$$\mathcal{H}(x,D\psi_1(x)) \le -\gamma, \qquad x \in \mathcal{W}. \tag{5.11}$$

For $\beta \in (0,1)$ define

$$\check{v}_\beta \doteq \beta\check{v} + (1-\beta)\psi_1. \tag{5.12}$$

It suffices to show that, for every $\beta \in (0,1)$,

$$\check{v}_\beta(x) \le v(x), \qquad x \in \mathcal{W}. \tag{5.13}$$

We argue by contradiction and assume that (5.13) does not hold. Therefore, there exist $\beta \in (0,1)$ and $\xi \in \mathcal{W}$ such that

$$\check{v}_\beta(\xi) - v(\xi) = \max_{x\in\mathcal{W}}(\check{v}_\beta(x) - v(x)) \doteq \delta > 0. \tag{5.14}$$

We will first argue that if $z \in \mathcal{W}^o$ and $(q,A) \in \overline{J}^{2,+}\check{v}_\beta(z)$, then

$$F(z,\check{v}_\beta(z),q,A) \le -(1-\beta), \qquad \mathcal{H}(x,q) \le -(1-\beta)\gamma. \tag{5.15}$$

To this end, consider first, $(q,A) \in J^{2,+}\check{v}_\beta(z)$. Let $\varphi \in C^2(\mathcal{W})$ be such that $q = D\varphi(z)$, $A = D^2\varphi(z)$ and $\check{v}_\beta - \varphi$ has a local maximum at $z$. Let

$$\widetilde{\varphi} = \beta^{-1}(\varphi - (1-\beta)\psi_1).$$

Then $\check{v} - \widetilde{\varphi} = \beta^{-1}(\check{v}_\beta - \varphi)$ has a local maximum at $z$. Let $q^*$ and $A^*$ be such that $q^* = D\widetilde{\varphi}(z)$, $A^* = D^2\widetilde{\varphi}(z)$. Note that

$$q = \beta q^* + (1-\beta)D\psi_1(z), \qquad A = \beta A^* + (1-\beta)D^2\psi_1(z). \tag{5.16}$$

Then $(q^*,A^*) \in J^{2,+}\check{v}(z)$. Using the subsolution property of $\check{v}$ (cf. Definition 2.2 and the text immediately following it),

$$F(z,\check{v}(z),q^*,A^*) \le 0. \tag{5.17}$$

Noting that the map $(r,u,X)\mapsto F(z,r,u,X)$ is affine and combining (5.10), (5.12), (5.16) and (5.17), we obtain the first inequality in (5.15). Now since $z\in\mathcal{W}^o$ and $(q,A)\in J^{2,+}\check{v}_\beta(z)$ are arbitrary, and $F$ is continuous in all variables, the first inequality in (5.15) holds, in fact, for all $z\in\mathcal{W}^o$ and $(q,A)\in\overline{J}^{2,+}\check{v}_\beta(z)$.

By convexity of $q\mapsto\mathcal{H}(z,q)$, (5.11) and using once more the subsolution property of $\check{v}$ and the continuity of $\mathcal{H}$, it is seen that

$$\mathcal{H}(z,q)\le\beta\mathcal{H}(z,q^*)+(1-\beta)\mathcal{H}(z,D\psi_1(z))\le-(1-\beta)\gamma,$$

for all $(q,A)\in\overline{J}^{2,+}\check{v}_\beta(z)$ and $z\in\mathcal{W}^o$. This proves the second inequality in (5.15).

Recall that $\xi$ is chosen such that (5.14) holds. Let $\eta=\eta(\xi)$ be as in (5.1). For $\widetilde{\gamma}\in(1,\infty)$ and $\varepsilon\in(0,1)$ set

$$\Psi(x,y)=|\widetilde{\gamma}(x-y)-\varepsilon\eta|^2+\varepsilon|y-\xi|^2,$$
$$\Phi(x,y)=\check{v}_\beta(x)-v(y)-\Psi(x,y),\qquad (x,y)\in\mathcal{W}\times\mathcal{W},$$

and let

$$(\widetilde{x}_{\varepsilon,\widetilde{\gamma}},\widetilde{y}_{\varepsilon,\widetilde{\gamma}})\equiv(\widetilde{x},\widetilde{y})\in\mathop{\arg\max}_{(x,y)\in\mathcal{W}\times\mathcal{W}}\Phi(x,y)$$
$$\equiv\left\{(u,w):\Phi(u,w)=\max_{(x,y)\in\mathcal{W}\times\mathcal{W}}\Phi(x,y)\right\}.$$

By (5.1),

$$\xi+\frac{\varepsilon}{\widetilde{\gamma}}\eta\in\mathcal{W}^o. \tag{5.18}$$

Clearly, $\Phi(\widetilde{x},\widetilde{y})\ge\Phi(\xi+\widetilde{\gamma}^{-1}\varepsilon\eta,\xi)$. This can be rewritten as

$$\check{v}_\beta(\widetilde{x})-v(\widetilde{y})-\check{v}_\beta\left(\xi+\frac{\varepsilon}{\widetilde{\gamma}}\eta\right)+v(\xi)\ge|\widetilde{\gamma}(\widetilde{x}-\widetilde{y})-\varepsilon\eta|^2+\varepsilon|\widetilde{y}-\xi|^2. \tag{5.19}$$

Dividing by $\widetilde{\gamma}^2$, we see that, for every $\varepsilon$, $|\widetilde{x}-\widetilde{y}|\to 0$ as $\widetilde{\gamma}\to\infty$. This observation, along with (5.14), (5.19) and the continuity of $\check{v}_\beta$ and $v$, gives that $\limsup_{\widetilde{\gamma}\to\infty}|\widetilde{\gamma}(\widetilde{x}-\widetilde{y})-\varepsilon\eta|^2+\varepsilon|\widetilde{y}-\xi|^2\le 0$. Hence, for all $\varepsilon\in(0,1)$,

$$\widetilde{y}\to\xi,\ \widetilde{\gamma}(\widetilde{x}-\widetilde{y})\to\varepsilon\eta\qquad\text{as }\widetilde{\gamma}\to\infty. \tag{5.20}$$

In particular,

$$\widetilde{x}=\widetilde{y}+\widetilde{\gamma}^{-1}\varepsilon\eta+\widetilde{\gamma}^{-1}o(1) \tag{5.21}$$

as $\widetilde{\gamma}\to\infty$. Hence, by (5.18) and (5.1), $\widetilde{x}\in\mathcal{W}^o$ for $\widetilde{\gamma}>\widetilde{\gamma}_0$, for some $\widetilde{\gamma}_0=\widetilde{\gamma}_0(\varepsilon)<\infty$. By (5.10), (5.12) and (5.14), it follows that $v(\xi)<\check{v}(\xi)$. By choosing $\widetilde{\gamma}_0$ larger if necessary, we have $v(\widetilde{y})<\check{v}(\widetilde{y})$ for $\widetilde{\gamma}>\widetilde{\gamma}_0$. Henceforth, assume $\widetilde{\gamma}>\widetilde{\gamma}_0$. For $(x,r,q,A)\in\mathcal{W}\times\mathbb{R}\times\mathbb{R}^d\times\mathcal{S}(d)$ let

$$\check{F}(x,r,q,A)\doteq F(x,r,q,A)\vee\mathcal{H}(x,q).$$

Since $\widetilde{x}\in\mathcal{W}^o$, we have from (5.15)

$$\check{F}(\widetilde{x},\check{v}_\beta(\widetilde{x}),q,X)\le -\varepsilon^*,\qquad (q,X)\in\overline{J}^{2,+}\check{v}_\beta(\widetilde{x}),$$

where $\varepsilon^*\doteq(1-\beta)(1\wedge\gamma)$. By the supersolution property of $v$,

$$\check{F}(\widetilde{y},v(\widetilde{y}),q,Y)\ge 0\qquad \text{for all } (q,Y)\in\overline{J}^{2,-}v(\widetilde{y}).$$

Combine the above two displays and use the inequality $(a\vee b)-(c\vee d)\le(a-c)\vee(b-d)$, $(a,b,c,d)\in\mathbb{R}^4$, along with (5.2)–(5.6), to obtain, for all $(q_1,X)\in\overline{J}^{2,+}\check{v}_\beta(\widetilde{x})$ and $(q_2,Y)\in\overline{J}^{2,-}v(\widetilde{y})$,

$$\varepsilon^*\le\check{F}(\widetilde{y},v(\widetilde{y}),q_2,Y)-\check{F}(\widetilde{x},\check{v}_\beta(\widetilde{x}),q_1,X)\le\Delta_1\vee\Delta_2, \tag{5.22}$$

where

$$\begin{aligned}\Delta_1&=\omega(|\widetilde{x}-\widetilde{y}|)+C|q_1-q_2|,\\ \Delta_2&=C^2|\widetilde{x}-\widetilde{y}|+\alpha(\widetilde{x})(v(\widetilde{y})-\check{v}_\beta(\widetilde{x}))+C|\widetilde{x}-\widetilde{y}|\,|q_1|\\ &\qquad+C|q_1-q_2|+\omega(|\widetilde{x}-\widetilde{y}|)\\ &\qquad+\tfrac{1}{2}\mathrm{trace}(\Gamma(\widetilde{x})X-\Gamma(\widetilde{y})Y).\end{aligned} \tag{5.23}$$

With an abuse of notation, we used the symbol $C$ in the above display for $C\vee\max_{\mathcal{W}}|v|$ (and will keep this notation). Next, noting that

$$\check{v}_\beta(\widetilde{x})-v(\widetilde{y})\ge\Phi(\widetilde{x},\widetilde{y})\ge\Phi(\xi,\xi),$$

and using (5.14), we have

$$v(\widetilde{y})-\check{v}_\beta(\widetilde{x})\le\varepsilon^2|\eta|^2. \tag{5.24}$$

Hence, by (5.23),

$$\begin{aligned}\Delta_2\le\ &C^2|\widetilde{x}-\widetilde{y}|(1+|q_1|)+C\varepsilon^2|\eta|^2+C|q_1-q_2|+\omega(|\widetilde{x}-\widetilde{y}|)\\ &+\tfrac{1}{2}\mathrm{trace}(\Gamma(\widetilde{x})X-\Gamma(\widetilde{y})Y).\end{aligned} \tag{5.25}$$

We now estimate the last term in the above display. By Theorem 3.2 of [3], since $\Phi$ has a (local) maximum at $(\widetilde{x},\widetilde{y})$, for each $\varrho\in(0,\infty)$, one can find $X,Y\in\mathcal{S}(d)$ such that

$$(D_x\Psi(\widetilde{x},\widetilde{y}),X)\in\overline{J}^{2,+}\check{v}_\beta(\widetilde{x}),\qquad(-D_y\Psi(\widetilde{x},\widetilde{y}),Y)\in\overline{J}^{2,-}v(\widetilde{y}),$$

and, with the usual order on $\mathcal{S}(2d)$,

$$\begin{pmatrix}X&0\\0&-Y\end{pmatrix}\le D^2\Psi(\widetilde{x},\widetilde{y})+\varrho(D^2\Psi(\widetilde{x},\widetilde{y}))^2. \tag{5.26}$$

Observing that

$$\begin{aligned}D_x\Psi(\widetilde{x},\widetilde{y})&=2\widetilde{\gamma}(\widetilde{\gamma}(\widetilde{x}-\widetilde{y})-\varepsilon\eta),\\ -D_y\Psi(\widetilde{x},\widetilde{y})&=2\widetilde{\gamma}(\widetilde{\gamma}(\widetilde{x}-\widetilde{y})-\varepsilon\eta)-2\varepsilon(\widetilde{y}-\xi)\end{aligned} \tag{5.27}$$

and

$$D^2\Psi(\tilde{x},\tilde{y}) = 2\tilde{\gamma}^2\begin{pmatrix} I & -I \\ -I & I \end{pmatrix} + 2\varepsilon\begin{pmatrix} 0 & 0 \\ 0 & I \end{pmatrix}, \tag{5.28}$$

we can rewrite (5.26) as

$$\begin{aligned}\begin{pmatrix} X & 0 \\ 0 & -Y \end{pmatrix} \le{}& (2\tilde{\gamma}^2 + 8\varrho\tilde{\gamma}^4)\begin{pmatrix} I & -I \\ -I & I \end{pmatrix} \\ &+ 4\tilde{\gamma}^2\varepsilon\varrho\begin{pmatrix} 0 & -I \\ -I & 2I \end{pmatrix} + (2\varepsilon + 4\varrho\varepsilon^2)\begin{pmatrix} 0 & 0 \\ 0 & I \end{pmatrix}.\end{aligned} \tag{5.29}$$

Note that if $A, B \in \mathcal{S}(2d)$ are nonnegative then $\operatorname{trace}(AB) \ge 0$. Arguing similarly to Example 3.6 [3], we use this fact along with (5.29) and the nonnegativity of the symmetric matrix

$$\begin{pmatrix} \sigma(\tilde{x})\sigma(\tilde{x})' & \sigma(\tilde{y})\sigma(\tilde{x})' \\ \sigma(\tilde{x})\sigma(\tilde{y})' & \sigma(\tilde{y})\sigma(\tilde{y})' \end{pmatrix}$$

to obtain

$$\begin{aligned}&\operatorname{trace}(\Gamma(\tilde{x})X - \Gamma(\tilde{y})Y) \\ &\quad = \operatorname{trace}(\sigma(\tilde{x})\sigma(\tilde{x})'X - \sigma(\tilde{y})\sigma(\tilde{y})'Y) \\ &\quad \le (2\tilde{\gamma}^2 + 8\varrho\tilde{\gamma}^4)\operatorname{trace}((\sigma(\tilde{x}) - \sigma(\tilde{y}))(\sigma(\tilde{x})' - \sigma(\tilde{y})')) \\ &\qquad + 8\tilde{\gamma}^2\varepsilon\varrho\operatorname{trace}(\sigma(\tilde{y})(\sigma(\tilde{y})' - \sigma(\tilde{x})')) + (2\varepsilon + 4\varrho\varepsilon^2)\operatorname{trace}(\sigma(\tilde{y})\sigma(\tilde{y})') \\ &\quad \le (2\tilde{\gamma}^2 + 8\varrho\tilde{\gamma}^4)\bar{C}|\tilde{x} - \tilde{y}|^2 + 8\bar{C}\tilde{\gamma}^2\varepsilon\varrho|\tilde{x} - \tilde{y}| + (2\varepsilon + 4\varrho\varepsilon^2)\bar{C},\end{aligned}$$

where $\bar{C} = (dC)^2$. By (5.22), (5.23), (5.25), (5.27) and the above estimate, substituting $\varrho = \tilde{\gamma}^{-2}$, we have for $\varepsilon < 1$, $\tilde{\gamma} > (1 \vee \tilde{\gamma}_0)$,

$$\begin{aligned}\varepsilon^* \le{}& [\bar{C}(1 + 2\tilde{\gamma}o(1)) + 4\bar{C}\varepsilon]|\tilde{x} - \tilde{y}| + 5\tilde{\gamma}^2\bar{C}|\tilde{x} - \tilde{y}|^2 \\ &+ 2\omega(|\tilde{x} - \tilde{y}|) + 4\bar{C}\varepsilon(1 + |\tilde{y} - \xi|) + \bar{C}\varepsilon^2|\eta|^2,\end{aligned}$$

where, as in (5.21), we wrote $o(1)$ for a function (possibly depending on $\varepsilon$) converging to zero as $\tilde{\gamma} \to \infty$. Let $\tilde{\gamma} \to \infty$ and use (5.20) and (5.21) to obtain

$$\varepsilon^* \le 6\bar{C}\varepsilon^2|\eta|^2 + 4\bar{C}\varepsilon.$$

Finally, letting $\varepsilon \to 0$ we arrive at a contradiction. Hence, (5.13) and, in turn, (5.9) must hold, and the result follows. □

## APPENDIX

LEMMA A.1. *The following are equivalent:*

(i) $\mathcal{H}(0) < \infty$.
(ii) $\mathcal{H}(q) < \infty$ *for some* $q \in \mathbb{R}^d$.
(iii) $\mathcal{H}(q) < \infty$ *for all* $q \in \mathbb{R}^d$.
(iv) *There exists* $c_1 \in (0, \infty)$ *such that* $(\kappa \cdot u)^- \le c_1|Gu|$ *for all* $u \in \mathcal{U}$.

(v) *The set* $\{u \in \mathcal{U} : |Gu| \le \varepsilon, \kappa \cdot u \le -1\}$ *is empty for all* $\varepsilon > 0$ *sufficiently small.*

PROOF. The equivalence of (i), (ii) and (iii) is immediate on noting that $|Gu \cdot q| \le |q|$ for all $u \in \mathcal{U}_1$.

Note that if for some $c \in (0,\infty)$ one has, for every $u \in \mathcal{U}$ with $Gu \neq 0$,

$$(\kappa \cdot u)^- \le c|Gu|, \tag{A.1}$$

then (A.1) holds for every $u \in \mathcal{U}$. Indeed, if $Gu = 0$, let $\bar{u}$ be a vector in $\mathcal{U}$ such that $|G\bar{u}| = 1$ [$\bar{u}$ exists by (1.2)]. Then for $r > 0$, $u^r \doteq u + r\bar{u} \in \mathcal{U}$ and $|Gu^r| > 0$, and so (A.1) holds for $u^r$ in place of $u$, and sending $r \to 0$, it follows that (A.1) holds for $u$ as well.

Note that the following holds: $\mathcal{H}(0) \le \sup_{u \in \mathcal{U}_1} (\kappa \cdot u)^- \le (\mathcal{H}(0))^+$. The implication (iv) $\Rightarrow$ (i) is immediate from the first inequality above. Conversely, if (i) holds, then the second inequality above, along with the argument of the last paragraph above, gives $(\kappa \cdot u)^- \le (\mathcal{H}(0))^+ |Gu|$ for every $u \in \mathcal{U}$, and (iv) follows.

Suppose now that (v) holds. Then there exists an $\varepsilon > 0$ such that for $u \in \mathcal{U}$ satisfying $|Gu| \le \varepsilon$ we have $(\kappa \cdot u)^- < 1$. In particular, for any $u \in \mathcal{U}$ with $|Gu| \neq 0$,

$$(\kappa \cdot u)^- \le \varepsilon^{-1} |Gu|. \tag{A.2}$$

Thus, by the argument given in the second paragraph of this proof, (iv) holds. Conversely, suppose that (iv) holds. Then $(\kappa \cdot u)^- < 1$ for all $u \in \mathcal{U}$ satisfying $|Gu| \le (2c_1)^{-1}$. In particular, the set in (v) with $\varepsilon = (2c_1)^{-1}$ is empty. We have thus established the equivalence of (i)–(v). □

LEMMA A.2. *Condition (1.8) implies (1.9).*

PROOF. Note that the result holds trivially if $\kappa = 0$. Assume now that $\kappa \neq 0$ and suppose that (1.9) fails. Then we can find a sequence $\{u_n\}$ contained in $\mathcal{U}$ such that $Gu_n \to 0$ and $\kappa \cdot u_n \le -1$. Note that $|u_n|$ is bounded from below by $|\kappa|^{-1}$. Let $\hat{u}_n = u_n |u_n|^{-1}$. Clearly, $G\hat{u}_n \to 0$ and $\kappa \cdot \hat{u}_n \le 0$. Let $\hat{u}_n$ converge to $\hat{u}$ along some subsequence. Then $\hat{u}$ is a unit vector satisfying $G\hat{u} = 0$ and $\kappa \cdot \hat{u} \le 0$. Thus, (1.8) fails. This proves the result. □

PROOF OF PART (ii) OF LEMMA 2.1. Note that by item (iv) of Lemma A.1, $(\kappa \cdot U_t)^- \le c_1 |GU_t| \ \le c_1 |A_t|$, for all $t \ge 0$, where

$$A_t = W_t - w_0 - \int_0^t \vartheta(W_s)\, ds - \int_0^t \sigma(W_s)\, dZ_s.$$

From the boundedness of $\mathcal{W}$ we have that, a.s., $\int_{[0,\infty)} e^{-\alpha t} |A_t|\, dt < \infty$, and $e^{-\alpha t} A_t \to 0$ as $t \to \infty$. Thus, we obtain that a.s., $\int_{[0,\infty)} \alpha e^{-\alpha s} (\kappa \cdot U_s)^-\, ds < \infty$ and $e^{-\alpha t}(\kappa\cdot U_t)^-\to 0$ as $t\to\infty$. In order to prove the assertion, it therefore suffices by (2.2) to show that if $\int_{[0,\infty)}\alpha e^{-\alpha t}(\kappa\cdot U_t)^+\,dt<\infty$, then $e^{-\alpha t}(\kappa\cdot U_t)^+\to 0$ as $t\to\infty$. The proof of this last statement follows the proof of Lemma A.4 of [7]. In particular, replacing $v$ by $\kappa$ and $Y$ by $U$ therein, we obtain in case (a) of the proof of that lemma [and the subcase where (100) there holds] that, for all $n$ sufficiently large,

$$\kappa\cdot\Delta_n U<\frac{-\varepsilon b}{1+\varepsilon}e^{\alpha\tau_{2n+1}}<0,$$

using the notation of [7]. Thus,

$$\begin{aligned}\frac{\varepsilon b}{1+\varepsilon}&\le(\kappa\cdot\Delta_n U)^- e^{-\alpha\tau_{2n+1}}\\&\le c_1|G\Delta_n U|e^{-\alpha\tau_{2n+1}}\\&\le c_2(|\Delta_n\check{Z}|+|\Delta_n W|+\Delta_n t)e^{-\alpha\tau_{2n+1}},\end{aligned}$$

where $\check{Z}_t\doteq\int_0^t\sigma(W_u)\,dZ_u$ and $c_2$ is a suitable positive constant. Since a.s., $\tau_{2n+1}\to\infty$ as $n\to\infty$, and $\Delta_n t=\tau_{2n+2}-\tau_{2n+1}<\delta$ for all but finitely many $n$, the right side above converges a.s. to 0, as $n\to\infty$, due to the compactness of the state space and the asymptotic properties of Brownian motion. Thus, $\frac{\varepsilon b}{1+\varepsilon}\le 0$, which is a contradiction since both $\varepsilon$ and $b$ are positive. Other cases are treated exactly as in [7], Lemma A.4. □

LEMMA A.3. (i) *Let $w\in\mathcal{W}$ and $u\in\mathcal{U}$. If $w+Gu\in\mathcal{W}$, then*

$$V(w+Gu)+\kappa\cdot u\ge V(w).$$

(ii) *Let $\tau$, $\varepsilon$, $w_0$ and $\mathcal{P}_{w_0,\varepsilon}$ be as in the proof of Proposition 4.2. Then $\inf_{\mathbb{P}^*\in\mathcal{P}_{w_0,\varepsilon}}\hat{J}(\mathbb{P}^*)=\inf_{\mathbb{P}^*\in\mathcal{P}_{w_0}}\hat{J}(\mathbb{P}^*)$, where $\hat{J}(\mathbb{P}^*)$ is as defined in (4.27).*

PROOF. (i) Let $U\in\mathcal{A}(w+Gu)$. Then $\overline{U}=u+U\in\mathcal{A}(w)$. Also,

$$V(w)\le J(w,\overline{U})\;=\;J(w+Gu,U)+\kappa\cdot u.$$

Since $U\in\mathcal{A}(w+Gu)$ is arbitrary, the result follows.

(ii) Fix $\mathbb{P}^*\in\mathcal{P}_{w_0}$. Define stochastic processes $(W,U)$ on $(\mathcal{E},\mathcal{B}(\mathcal{E}))$ as

$$W(s)\doteq\pi_1(s),\qquad U(s)\doteq\pi_2(s),\qquad s\in[0,\tau),$$

and on the set $\{\tau<\infty\}$ define

$$\varepsilon^*=\inf\{\delta\in[0,1]:W(\tau-)+\delta G\Delta\pi_2(\tau)\notin B_{2\varepsilon}(w_0)\}$$

when this latter set is nonempty, otherwise set $\varepsilon^*=1$. Let

$$W(\tau)=W(\tau-)+\varepsilon^*G\Delta\pi_2(\tau),\qquad U(\tau)\doteq\pi_2(\tau-)+\varepsilon^*\Delta\pi_2(\tau).$$

We leave $(W,U)$ unspecified on $(\tau,\infty)$; they can be defined in an arbitrary way, as long as $U \in \mathcal{A}(w_0)$ with system $\Phi^* \doteq (\mathcal{E},\mathcal{B}(\mathcal{E}),\{\mathcal{F}_t\},\mathbb{P}^*,\pi_3)$, $\tau = \inf\{t \ge 0 : W(t) \notin \overline{B_\varepsilon(w_0)}\}$, and (4.2) holds. Let $\tilde{\mathbb{P}}^*$ be the measure induced by $(W,U,\pi_3)$ on $(\mathcal{E},\mathcal{B}(\mathcal{E}))$. Then $\tilde{\mathbb{P}}^* \in \mathcal{P}_{w_0,\varepsilon}$. Finally, note that, setting $\rho \doteq (1-\varepsilon^*)\Delta\pi_2(\tau)1_{\{\tau<\infty\}}$, we have from the first part of the lemma that

$$\hat{J}(\mathbb{P}^*) - \hat{J}(\tilde{\mathbb{P}}^*) = \mathbb{E}^*[1_{\{\tau \le t\}}\, e^{-\alpha\tau}[\kappa \cdot \rho + V(\pi_1(\tau)) - V(\pi_1(\tau) - G\rho)]] \ge 0.$$

This proves the lemma. □

The following well-known result is included for the sake of completeness. For a proof, we refer the reader to Theorem 10.2.2 (page 345) of [4].

LEMMA A.4. *Let $(\Omega,\mathcal{F},P)$ be a probability space, $\mathcal{T}$ be a Polish space, $\mathcal{B}(\mathcal{T})$ be the Borel sigma field on $\mathcal{T}$, and $X : \Omega \to \mathcal{T}$ be a measurable function. Let $\mathcal{G}$ be a sub sigma field of $\mathcal{F}$. Then a regular conditional probability distribution for $X$ given $\mathcal{G}$ exists, that is, there is a function $P : \Omega \times \mathcal{B}(\mathcal{T}) \to [0,1]$ such that, for $P$-a.e. $\omega \in \Omega$, $P(\omega,\cdot)$ is a probability measure on $(\mathcal{T},\mathcal{B}(\mathcal{T}))$; for each $A \in \mathcal{B}(\mathcal{T})$, $P(\cdot,A)$ is $\mathcal{G}$-measurable, and for all $C \in \mathcal{G}$ and $A \in \mathcal{B}(\mathcal{T})$, $P(C \cap \{X \in A\}) = \int_C P(\omega,A)\,dP(\omega)$.*

R. ATAR
DEPARTMENT OF ELECTRICAL ENGINEERING
TECHNION–ISRAEL INSTITUTE OF TECHNOLOGY
HAIFA 32000
ISRAEL
E-MAIL: atar@ee.technion.ac.il

A. BUDHIRAJA
DEPARTMENT OF STATISTICS
AND OPERATIONS RESEARCH
UNIVERSITY OF NORTH CAROLINA
CHAPEL HILL, NORTH CAROLINA 27599-3260
USA
E-MAIL: budhiraj@email.unc.edu

R. J. Williams
Department of Mathematics
University of California, San Diego
La Jolla, California 92093-0112
USA
E-mail: williams@math.ucsd.edu